\documentclass[11pt,a4paper,reqno]{amsart}
\usepackage{graphicx} %

\usepackage{times} %
\usepackage{amsmath} %
\usepackage{amssymb}  %
\usepackage{amsthm}
\usepackage{latexsym}
\usepackage{amsfonts,bbm}
\usepackage{xcolor}
\usepackage{mathtools}
\usepackage{enumerate}
\usepackage{cite}
\usepackage{tikz}
\usepackage{graphicx,caption}
\usepackage{todonotes}
\usepackage{float}
\usepackage{nicefrac}
\usepackage[foot]{amsaddr}
\usepackage{color}
\usepackage{graphicx}
\newcommand{\bE}{\mathbb E}

\usepackage{tikz}
\usetikzlibrary{shapes.geometric}
\usetikzlibrary{arrows.meta,arrows}

\usepackage{tabularx}
\usepackage{comment}
\usepackage[hidelinks, colorlinks=false]{hyperref}
\usepackage{cleveref}

\usepackage[normalem]{ulem}
\usepackage{algorithm}

\allowdisplaybreaks
\newtheorem{theorem}{Theorem}[section]
\newtheorem{remark}[theorem]{Remark}
\newtheorem{proposition}[theorem]{Proposition}
\newtheorem{corollary}[theorem]{Corollary}
\newtheorem{lemma}[theorem]{Lemma}

\newtheorem{definition}[theorem]{Definition}
\newtheorem{example}[theorem]{Example}

\newcommand{\calK}{\mathcal K}  
\newcommand{\calX}{\mathcal X}  
\newcommand{\calZ}{\mathcal Z}  
\newcommand{\calN}{\mathcal N}

\newcommand{\bbH}{\mathbb H}
\newcommand{\bbU}{\mathbb U}

\newcommand{\bK}{\mathbf K}
\newcommand{\R}{\ensuremath{\mathbb R}}  
\newcommand{\N}{\ensuremath{\mathbb N}}  

\renewcommand{\k}{\mathsf{k}}
\newcommand{\bk}{\mathbf k}

\newcommand{\linspan}{\operatorname{span}}
\newcommand{\dist}{\operatorname{dist}}
\DeclareMathOperator*{\argmin}{arg\,min}

\newcommand{\ceil}[1]{\lceil #1 \rceil}
\newcommand{\floor}[1]{\left\lfloor #1 \right\rfloor}

\usepackage{amsmath,amssymb,amsfonts}
\usepackage{latexsym}
\usepackage{mathrsfs}
\usepackage{bbm}
\usepackage{xcolor}
\usepackage{tikz}
\usepackage{url}
\usepackage{stmaryrd}
\usepackage{dsfont}
\usepackage{placeins}
\usepackage{lipsum}
\usepackage{amsfonts}
\usepackage{graphicx}
\usepackage{epstopdf}
\usepackage{amsopn}
\ifpdf
  \DeclareGraphicsExtensions{.eps,.pdf,.png,.jpg}
\else
  \DeclareGraphicsExtensions{.eps}
\fi

\newcommand{\braces}[1]{{\rm (}#1{\rm )}}

\newcommand{\wt}{\widetilde}
\newcommand{\wh}{\widehat}

\newcommand{\la}{\lambda}
\newcommand{\veps}{\varepsilon}

\allowdisplaybreaks

\title[Kernel-based Koopman approximants for control]{Kernel-based Koopman approximants for control: Flexible sampling, error analysis, and stability}
\author{Lea Bold$^{1}$} \address{$^{1}$Optimization-based Control Group, Institute of Mathematics, Technische Universität Ilmenau, Germany.\\ Mail: \textsc{\{lea.bold,friedrich.philipp, karl.worthmann\}@tu-ilmenau.de}}
\author{Friedrich M.\ Philipp$^{1}$}
\author{Manuel Schaller$^{2}$}\address{$^{2}$Faculty of Mathematics, Chemnitz University of Technology\\ Mail: \textsc{manuel.schaller@math.tu-chemnitz.de}}
\author{Karl Worthmann$^{1}$}
\thanks{F.M.P was funded by the Carl Zeiss Foundation within the project {\it DeepTurb--Deep Learning in and from Turbulence}. He was further supported by the free state of Thuringia and the German Federal Ministry of Education and Research (BMBF) within the project {\it THInKI--Th\"uringer Hochschulinitiative für KI im Studium}. L.B.\ gratefully acknowledges funding by the German Research Foundation (DFG) – Project-ID 545246093.}

\begin{document}

\begin{abstract}
    Data-driven techniques for analysis, modeling, and control of complex dynamical systems are on the uptake. Koopman theory provides the theoretical foundation for the popular kernel extended dynamic mode decomposition (kEDMD). In this work, we propose a novel kEDMD scheme to approximate nonlinear control systems accompanied by an in-depth error analysis. 
    {Key features are regularization-based robustness and an adroit decomposition into micro and macro grids enabling flexible sampling. But foremost}, we prove proportionality, i.e., explicit dependence on the distance to the (controlled) equilibrium, of the derived {bound} on the full approximation error. Leveraging this key property, we rigorously show that asymptotic stability of the data-driven surrogate (control) system implies asymptotic stability of the original (control) system and vice versa.\smallskip
    
\end{abstract}
\maketitle

\section{Introduction}

Data-driven methods for analysis, modelling, and control of dynamical systems {have} recently attracted considerable attention, see, e.g., the survey articles~\cite{prag2022toward,MartScho23} and the references therein. 
Key aspects are the Lyapunov-based stability analysis {of and the controller design for} dynamical systems, see, e.g., \cite{BerbKohl20}.
For nonlinear systems, however, guarantees for data-driven approaches require {more sophisticated techniques}~\cite{scharnhorst2022robust}. 
In this work, we provide a framework that allows to infer stability properties of the original system from its data-driven surrogates (and vice versa). 
Hereby, we leverage the Koopman operator, which provides a theoretically-sound foundation for analysis and control of dynamical systems through the lens of observable {functions}, see, e.g., \cite{rowley2009spectral,mauroy2016global,otto2021koopman}. 
{The underlying idea is to lift the dynamics to an infinite-dimensional function space. %
Then, %
finite-dimensional surrogate models are constructed using regression in a data-driven fashion}. 
The most prominent approximation scheme %
is extended dynamic mode decomposition (EDMD; \cite{WillKevr15}), an algorithm that builds upon a finite dictionary of observables and has been used in a wide range of challenging applications, see, e.g., \cite{XuWang23,ZhouLu24} for some recent examples.

Various works combining the Koopman approach with Lyapunov arguments exist in the literature. 
For a robust $\mathcal{H}_\infty$-approach, we refer to~\cite[Chapter 2]{MaurSusu20}, while Lyapunov-based stabilization of a chemical reactor was the subject of~\cite{NaraKwon19}. 
In~\cite{ZinaBako23}, a neural-network-based Koopman approach to feedback design via control Lyapunov functions is presented, while Lyapunov functions were approximated in a Koopman-based manner in~\cite{BreiHoev23}.
However, to rigorously ensure end-to-end guarantees, an error analysis including, e.g., convergence rates, of data-driven surrogate models serving as approximants is indispensable. 
Concerning EDMD, {strong} convergence to the Koopman operator in the infinite-data limit was proven in \cite{KordaMezi18} (using also infinitely-many observables). 
Finite-data error bounds for dynamical systems were, to the best of the authors' knowledge, first provided in~\cite{Mezi22,NuskPeit23,ZhanZuaz23}. 
Then, the probabilistic error bounds~\cite{Mezi22,ZhanZuaz23} were further extended to stochastic (control) systems using both i.i.d.\ and ergodic sampling in~\cite{NuskPeit23} and to kernel EDMD in~\cite{PhilScha24}, see also~\cite{kostic2024sharp} for sharp rates in view of spectral approximations and~\cite{kostic2024consistent} for long-term ergodic predictions via transfer operators.
Recently, using kernel EDMD (kEDMD) embedded in suitably-chosen reproducing kernel Hilbert spaces (RKHSs), error bounds in the supremum norm were rigorously shown in~\cite{KohnPhil24}. Similar estimates have been shown in~\cite{YadaMaur24} via an EDMD-variant using Bernstein polynomials.
A key finding in the derivation of these uniform error bounds is the invariance of the RKHS under the Koopman operator, see also~\cite{GonzAbud23} and~\cite{PoweParu24} for preliminary results.
Assuming such uniform error bounds, stabilizing controllers with end-to-end guarantees can be designed, see, e.g., \cite{GoswPale21,maddalena2021kpc,StraScha24:generator}, or \cite{BoldGrun23,WortStra24} within model predictive control. 
Herein, error bounds on EDMD-based surrogates for control systems that are proportional to the distance to the set point play a major role. 
While pointwise and proportional error bounds are fundamental in these works for controller design, there is, to the best of the authors' knowledge, no rigorous proof of such bound for an existing EDMD variant for control systems. 
This work closes this gap.
In this context, we point out~\cite{IacoToth24}, where limitations of linear surrogate models~\cite{KordMezi18:MPC, ProcBrun16} were thoroughly discussed.
Hence, for general nonlinear control-affine systems, the bilinear approach proposed in~\cite{Sura16,WillHema16}, is superior, see, e.g., \cite{PeitOtto20,GoswPale21}. 

The contribution of this work is two-fold. 
First, we propose a kEDMD-based surrogate model for nonlinear autonomous dynamics and prove that it inherits stability properties of the ground-truth system (and vice-versa). In particular, we show that asymptotic stability of one system implies practical asymptotic stability of the other. 
Further, we prove the first {uniform proportional} error bounds on kernel EDMD approximants of the Koopman operator, which certify that the full approximation error decays proportionally w.r.t.\ the distance to the equilibrium without imposing restrictive assumptions like invariance of some finitely-generated subspace~\cite{GoswPale21}. 
If a certain compatibility condition linking the decay rates of the stability-certifying Lyapunov function and the proportional error bound holds, we even prove that asymptotic stability is fully preserved.
As a second major contribution, we present novel data-based Koopman approximants using kEDMD, which allow for flexible sampling of state-control pairs and yield control-affine surrogate systems. 
We analyze our algorithm and prove proportional bounds on the uniform approximation error between the original model and its surrogate. 
The proportionality is then leveraged to prove that feedback laws stabilizing the data-driven model also stabilize the original control system (and vice versa). 
All results are extended in view of regularized kEDMD approximants to improve robustness, e.g., for noisy data.

The outline of the paper is as follows: 
In \Cref{sec:Koopman}, we begin with a brief recap on Koopman approximants for autonomous systems generated by kernel EDMD. 
Then, we present uniform bounds on the full approximation error using RKHS interpolation, before we conclude the section by proposing a data-driven surrogate and the respective error bound.
In the subsequent~\Cref{sec:stability}, we show that (practical) asymptotic stability is inherited 
from the original dynamical system by the data-driven surrogate 
and vice versa. 
In~\Cref{sec:kEDMD:control}, we propose a novel and highly-flexible approximation scheme for control systems, prove error bounds on the full approximation error and illustrate the results numerically.
Conclusions are drawn in~\Cref{sec:conclusions}.

\smallskip\noindent
{\bf Notation.} We use $\|\cdot\|$ for the Euclidean norm on $\R^n$ and its induced matrix norm on $\R^{n\times n}$. The Frobenius norm on $\R^{n\times m}$ will be denoted by $\|\cdot\|_F$. Moreover, we let $\R_{\geq 0} = [0,\infty)$. Further, {we set $\N = \mathbb Z_{\ge 1}$, and} for $d\in\N$ we use the abbreviation $[1:d] := \mathbb{Z} \cap [1,d]$. By $C_b(\Omega)$, we denote the space of bounded continuous functions on a set $\Omega\subset\R^n$. {\em Moduli of continuity} of a continuous function $f : \Omega\to\R$ are denoted by $\omega_f$. Recall that a modulus of continuity $\omega_f : \R_{\ge 0}\to\R_{\ge 0}$ vanishes at zero, is continuous at zero, and satisfies $| f(x) - f(y) | \leq \omega_f(\| x - y \|) \qquad\forall\, x,y \in \Omega$.
We define comparison functions to introduce our stability notions analogously to~\cite{GrunPann17}, see, e.g., \cite{Kell14}. 
A continuous, strictly increasing function $\alpha: \R_{\ge 0} \to \R_{\ge 0}$ satisfying $\alpha(0) = 0$ is said to be of class~$\mathscr{K}$. If it is, in addition, unbounded, the comparison function~$\alpha$ is class~$\mathscr{K}_\infty$. 
A continuous function $\beta: \R_{\ge 0}^2 \rightarrow \R_{\ge 0}$ is called a class~$\mathscr{K\!L}$-function if, for each $t \geq 0$, $\beta(\cdot,t) \in \mathscr{K}_\infty$ %
and $\beta(r,\cdot)$ is strictly monotonically decreasing with $\lim_{t \rightarrow \infty} \beta(r,t) = 0$ for all $r > 0$.

\section{Koopman operator and kernel EDMD}
\label{sec:Koopman}

Throughout this paper, let $\Omega \subset \mathbb{R}^n$ be an open domain with Lipschitz boundary in the sense of \cite[§4.9]{Adams2003}.\footnote{The imposed Lipschitz condition implies the usual cone conditions for interpolation estimates with reproducing kernel Hilbert spaces, cf.\ \cite[Appendix A]{KohnPhil24}.} We consider the discrete-time dynamical system given by 
\begin{align}\label{eq:dynamics}\tag{DS}
    x^+ = F(x)
\end{align}
with a map $F:\Omega\to\mathbb{R}^n$. 
Here, $x \in \Omega$ and $x^+$ are the current and the successor state of the dynamical system~\cref{eq:dynamics}, respectively.
Further, we abbreviate the image of the set~$\Omega$ w.r.t.\ the dynamics~$F$ by $\Xi$, i.e., $\Xi = F(\Omega) := \{F(x) : x\in\Omega\}$, and assume that $F \in \mathcal{C}^1(\Omega,\Xi)$ is a diffeomorphism satisfying the regularity condition $\inf_{x \in \Omega} | \det DF(x) | > 0$. 

\begin{remark}[Differential equation]\label{rem:ode}
    The dynamical system~\cref{eq:dynamics} can be inferred from the ordinary differential equation %
    \begin{align}\label{eq:ode}\tag{ODE}
        \dot{x}(t) = g_0(x(t)),
    \end{align}
    where $x(t) \in \Omega$ represents the state at time $t \geq 0$ and $g_0: \Omega \to \mathbb{R}^n$ is a locally Lipschitz-continuous map. 
    Then, for a given time step $\Delta t > 0$, we can associate the discrete-time dynamical system $x^+ = F(\bar{x}) := \bar{x} + \int_0^{\Delta t} g_0(x(s;\bar{x}))\,\mathrm{d}s$, using the integral representation of the solution $x(\cdot,\bar{x})$ emanating from the (current) state~$\bar{x} \in \Omega$ tacitly assuming that the solution $x(s;\bar{x})$ exists on the time interval $[0,\Delta t]$. 
    We point out that the imposed regularity assumption $\inf_{x \in \Omega} | \det DF(x) | > 0$ automatically holds for a sufficiently small time step~$\Delta t$.
\end{remark}

The linear Koopman operator~$\mathcal{K} : C_b(\Xi) \to C_b(\Omega)$ associated with the system~\cref{eq:dynamics} maps functions (so-called {\em observables}) in $C_b(\Xi)$ to functions in $C_b(\Omega)$ along the flow of the system and is defined by
$$
    (\mathcal{K} f)(x) = f(F(x)) \qquad\forall\,f\in C_b(\Xi),\;x\in\Omega,
$$
or, for short, $\mathcal{K} f = f\circ F$. 
{In the following, we assume forward invariance of~$\Omega$ w.r.t.\ the dynamics~\cref{eq:dynamics} in order to keep the presentation technically simple. In view of the stability analysis of the subsequent sections, we point out that this seemingly restrictive assumption holds on a suitably chosen sublevel set of the considered Lyapunov function.  For a detailed treatment of the case $\Xi = F(\Omega) \not\subset \Omega$, we refer to~\cite{KohnPhil24}.}

\subsection{Kernel EDMD (kEDMD)}\label{sec:sub:kEDMD}

Let $\k: \mathbb{R}^n \times \mathbb{R}^n \to \mathbb{R}$ be a continuous strictly positive-definite symmetric kernel, i.e., for any set of pairwise distinct data points
\begin{align}\label{eq:data-points}
    \calX = \{x_1,\ldots,x_d\} \subset \mathbb{R}^n,
\end{align}
the symmetric {\em kernel matrix} $(\k(x_i,x_j))_{i,j=1}^d$ is positive definite.\footnote{In the literature, a positive definite kernel is defined by positive semi-definiteness of the kernel matrices. Since we require these matrices to be invertible, we added the term {\em strictly}.} 
The {\it canonical feature}~$\phi_z$ of $\k$ at $z \in \mathbb{R}^n$ is defined by $\phi_z(x) = \k(z,x)$, $x \in \mathbb{R}^n$. 
It is well known that the linear space $\linspan\{ \phi_z : z \in \mathbb{R}^n \}$ extends by completion to a Hilbert space~$\bbH$ of functions, the {\em reproducing kernel Hilbert space} (RKHS) or simply {\em native space} associated with (or generated by) the kernel~$\k$, see, e.g., \cite[Theorem 2.14]{PaulRagh16}. Continuity of the kernel is inherited by the functions in~$\bbH$ so that $f|_\Omega\in C_b(\Omega)$ holds for all $f \in \bbH$. 
For all functions $f \in \bbH$, we have the {\em reproducing property}
\begin{align}\label{eq:RKHS:reproducing-property}
    f(x) = \langle f,\phi_x\rangle \qquad\forall\,x \in \mathbb{R}^n,
\end{align}
{where $\langle \cdot, \cdot \rangle$ is the inner product of the RKHS~$\mathbb{H}$. Before we proceed}, we provide an example of popular radially symmetric kernels with compact support.
\begin{example}[RKHS~$\mathbb{H}$ generated by Wendland kernels: $\mathcal{N}_{\Phi_{n,k}}$]\label{ex:wend}
    The Wendland radial basis function (RBF) $\Phi_{n,k} : \R^{n} \to \R$ of smoothness degree $k\in\N_0$ is defined by
    \[
        \Phi_{n,k}(x) := \phi_{n,k}(\|x\|) \qquad\text{ with }\qquad \phi_{n,k}(r) = \begin{cases}
            p_{n,k}(r), & 0 \le r \le 1, \\
            0, & 1 < r,
        \end{cases}
    \]
    where $p_{n,k}$ is {a univariate} polynomial of degree $\floor{\frac{n}{2}} + 3k + 1$ and $\phi_{n,k} \in \mathcal{C}^{2k}([0,\infty),\mathbb{R})$, see \cite[Theorem 9.13]{Wend04}. 
    The Wendland RBF~$\Phi_{n,k}$ induces the kernel $\k_{n,k}$ given by $\k_{n,k}(x,y) = \Phi_{n,k}(\|x-y\|)$ for $x,y \in \R^n$. 
    By $\calN_{\Phi_{n,k}}(\mathcal{O})$, we denote the native space corresponding to the RKHS~$\mathbb{H}$ generated by the Wendland kernel on any bounded open domain $\mathcal{O} \subset \R^{n}$ with Lipschitz boundary. 
    Here, denoting by $H^s(\mathcal{O})$ the $L^2$-Sobolev space of regularity order~$s$ on~$\mathcal{O}$, we have the identity
    \begin{align}\label{eq:identity:RKHS:Sobolev:Wendland}
        \mathcal{N}_{\Phi_{n,k}}(\mathcal{O}) = H^{\sigma_{n,k}}(\mathcal{O}) \qquad\text{ with }\qquad \sigma_{n,k} := \frac{n+1}{2} + k
    \end{align}
    with equivalent norms, see e.g.~\cite[Corollary 10.48]{Wend04} for integer Sobolev orders or \cite[Theorem 4.1]{KohnPhil24} for fractional Sobolev orders orders. 
\end{example}

To obtain error estimates with finitely many data points, in the following we assume that the set $\Omega$ is bounded. 
{Moreover, we choose the kernel~$\k$ such that the corresponding RKHS~$\mathbb{H}$ is invariant w.r.t.\ the action of the Koopman operator $\mathcal{K}: C_b(\Omega) \to C_b(\Omega)$, i.e., $\calK \,\bbH \subset \bbH$, a key property to derive bounds on the approximation error for kernel-based approximations of the Koopman operator. 
A suitable choice are Wendland kernels (see \cref{ex:wend}) as rigorously shown} in~\cite[Section~4.2]{KohnPhil24}. 
The authors of~\cite{KohnPhil24} leverage that these native spaces coincide with fractional Sobolev spaces with equivalent norms, cf.\ identity~\cref{eq:identity:RKHS:Sobolev:Wendland}. We note that the same holds true for {\em Matérn kernels}, as illustrated in~\cref{rem:Matern}.
{We point out that, in general, the invariance $\calK \,\bbH \subset \bbH$ does not hold for nonlinear systems, e.g., the RKHS~$\mathbb{H}$ generated by Gaussian kernels only satisfies this invariance condition if the dynamics~\eqref{eq:dynamics} are affine linear as shown in~\cite{GonzAbud23}}. 

Before introducing kernel-based approximations of the Koopman operator, we briefly note an important consequence of the reproducing property~\cref{eq:RKHS:reproducing-property} for the Koopman operator: 
A point $x^*\in \Omega$ is an equilibrium of the dynamics~\cref{eq:dynamics} if and only if $(\mathcal{K} f)(x^*) = f(x^*)$ for all $f\in\bbH$. Indeed, {necessity} is trivial by definition of the Koopman operator. For {sufficiency}, assume that $(\mathcal{K} f)(x^*) = f(x^*)$ for all $f\in\bbH$.
Then, for all $f\in\bbH$ we have $\langle f,\phi_{F(x^*)}\rangle  = f(F(x^*)) = (\mathcal{K} f)(x^*) = f(x^*) = \langle f,\phi_{x^*}\rangle$ and, thus, $\phi_{F(x^*)} = \phi_{x^*}$. The strict positive definiteness of the kernel~$\k$ then implies $F(x^*) = x^*$, i.e., that $x^*$ is an equilibrium of \cref{eq:dynamics}.

Next, we briefly recap {\em kernel extended dynamic mode decomposition} (kEDMD) as an advanced tool to approximate the Koopman operator from data and refer to~\cite{WillRowl15,KlusNusk20} and~\cite[Section 3.2]{KohnPhil24} for further details. For the set of pairwise distinct data points~$\mathcal{X}$ given by~\cref{eq:data-points}, we set $V_\calX = \operatorname{span} \{\phi_{x_1},\ldots,\phi_{x_d}\}$, where $\phi_{x_i} = \k(x_i,\cdot)$, $i \in [1:d]$. 
Further, let $P_\calX$ denote the orthogonal projection in $\bbH$ onto $V_\calX$, i.e., for given $f \in \bbH$, the function~$P_\calX f$ solves the regression problem
\begin{align}\label{eq:regression}
    \min_{g\in V_\calX} \|f-g\|_{\bbH}^2,
\end{align}
cf.~\cite{KohnPhil24}. Hence, $g^*\in V_\calX$ is the solution of~\cref{eq:regression} if and only if $g^* = P_\calX f$. 
Further, it may be easily seen that \cref{eq:regression} is equivalent to $\min_{g\in V_\calX} \sum_{i=1}^d |f(x_i) - g(x_i)|^2$. Note that $P_\calX f$ for $f \in \bbH$ is in fact the interpolation of~$f$ at the points in~$\calX$ in the sense that $P_\calX f$ is the unique function in $V_\calX$ which coincides with $f$ at all data points in~$\calX$.
Then, as proven in \cite[Proposition~3.2]{KohnPhil24}, a kEDMD approximant of the Koopman operator~$\mathcal{K}$ on~$\bbH$ is given by
\begin{align}\label{eq:edmd}
    \widehat{\mathcal{K}} = P_{\mathcal{X}} \mathcal{K}.
\end{align}
Now, letting $f_\calZ := [f(z_1),\ldots,f(z_d)]^\top$ for a set $\calZ = \{z_1,\ldots,z_d\}$, the approximant $\widehat{\mathcal{K}} f$
of $\mathcal{K} f$ may be written as
\begin{align}\label{e:matrices}
\widehat{\mathcal{K}} f = P_\calX \mathcal{K} f = P_\calX (f\circ F) = f_{F(\calX)}^\top\bK_\calX^{-1}\bk_\calX,
\end{align}
where $\bk_\calX = [\phi_{x_1},\ldots,\phi_{x_d}]^\top$ and $\bK_\calX = (k(x_i,x_j))_{i,j=1}^d$.
In fact, the kEDMD regression solution is $P_\calX \mathcal{K}|_{V_\calX}$, which is a linear map from the finite-dimensional space $V_\calX$ into itself. 
Correspondingly, kEDMD may also be understood as a method for approximating the Koopman operator by finite-rank operators with range in $V_\calX$ despite the Koopman operator acting on an infinite-dimensional function space. 

We close this subsection by relating the above approximation~\cref{eq:edmd} to another kernel-based surrogate of the Koopman operator. 
\begin{remark}\label{rem:alternative}
    In~\cite{KohnPhil24}, another approximant of $\calK$ has been defined by $\widetilde \calK = P_\calX\calK P_\calX$. 
    The key difference is the following: 
    For the computation of $\widehat{\calK} f$, the observable $f\in\bbH$ has to be propagated by the flow, i.e., we require data samples $(f(x_i),f(F(x_i))$, $i \in [1:d]$.
    In contrast, for the computation of the alternative surrogate~$\widetilde{\calK} f$, the canonical features $\phi_{x_i}$ at the data sites have to be propagated, while the observable~$f$ only has to be interpolated, i.e., measurements $f(F(x_i))$ are not necessary. 
    We refer the interested reader to~\cite{KohnPhil24} for a detailed discussion. %
    Moreover, a collection of results on the surrogate arising from the approximant~$\widetilde\calK$ --~which correspond to those that we will derive for~$\widehat\calK$~-- can be found in the (additional) appendix, see \Cref{s:K_tilde}.
\end{remark}

\subsection{Approximation error for regularized kEDMD: uniform bounds}
\label{sec:sub:error:bound}

In this section, we extend recently proposed uniform error bounds for the approximant $\widehat\calK$ in the operator norm \mbox{$\|\cdot\|_{\bbH\to C_b(\Omega)}$} to regularized kEDMD. 

If the size of the data set $\calX$ is large, the evaluation of \cref{e:matrices} may lead to numerical instabilities as the kernel matrix $\bK_\calX$ is typically badly conditioned. 
For this reason, one often regularizes %
kernel-based interpolation problems of the form \cref{eq:regression}. More precisely, the regression problem is endowed by a regularization term, i.e., 
\begin{align*}
    R^\lambda_\calX f := \argmin\nolimits_{g\in V_\calX}\;\sum\nolimits_{i=1}^d|f(x_i) - g(x_i)|^2 + \lambda\|g\|_{\calN_{\Phi_{n,k}}(\Omega)}^2
\end{align*}
with regularization parameter $\lambda\geq 0$. Clearly, $R_\calX^0 = P_\calX$ holds and it can be verified that the solution operator $R_\calX^\lambda$ is linear and satisfies $R^\lambda_\calX f = f_\calX^\top(\bK_\calX + \lambda I)^{-1}\bk_\calX$, $f\in\calN_{\Phi_{n,k}}(\Omega)$. 
However, while $P_\calX$ is a projection, the operator $R^\lambda_\calX$, $\lambda > 0$, is not. Nevertheless, as proved in the (additional) appendix (see \Cref{{p:RX}}), $R_\calX^\lambda$, $\lambda \geq 0$, shares with $P_\calX$ that it is self-adjoint and positive semi-definite as an operator on $\calN_{\Phi_{n,k}}(\Omega)$. Moreover, $R_\calX^\lambda$, $\lambda \geq 0$, commutes with $P_\calX$.
For regularized kEDMD, we define the following approximant of the Koopman operator: 
\begin{align}\label{eq:koopman_reg}
    \widehat{\mathcal{K}}_\lambda f = f_{F(\calX)}^\top(\bK_\calX + \lambda I)^{-1}\bk_\calX,
\end{align}
where $\lambda \geq 0$. It is not hard to see that $\widehat\calK_\la = R_\calX^\lambda\calK$. Moreover, the approximation~\cref{eq:edmd} is recovered for $\lambda = 0$, that is, $\widehat{\mathcal{K}}_0 = \widehat{\mathcal{K}}$. 

The next theorem provides a novel bound on the approximation error. This result extends the previous work~\cite{KohnPhil24}, which considered the case $\lambda = 0$, i.e., non-regularized interpolation. For a data set~$\calX\subset \Omega$ as in~\cref{eq:data-points}, we denote the {\em fill distance} by
\begin{align*}
    h_{\calX} := \sup\nolimits_{x \in \Omega}\,\dist(x,\calX)
\end{align*}
using the Euclidean norm $\|\cdot\|$ as the metric for the distance operator. {Analogously to~\cite{KohnPhil24}, and refer to~\cref{rem:Matern} for a discussion on a potential alternative.}
\begin{theorem}\label{thm:error:regularization}
    Let $k\geq 1$, $\lambda \geq 0$ and $F\in C_b^{\ceil{\sigma_{n,k}}}(\Omega;\R^{n})$. Then there are constants $C,h_0>0$ such that for any finite set $\calX = \{x_i\}_{i=1}^d \subset\Omega$ of {pairwise distinct} sample points with $h_{\calX}\le h_0$ and for all {$f\in\calN_{\Phi_{n,k}}(\Omega)$ and $x\in\Omega$}, we have
    \begin{align*}    
        |(\mathcal{K}f)(x) - (\widehat{\mathcal{K}}_\lambda f)(x)| \leq C\big(h_\calX^{k+1/2} + \sqrt\lambda\big)\|f\|_{\calN_{\Phi_{n,k}}(\Omega)},
    \end{align*}
    {where the RKHS $\calN_{\Phi_{n,k}}(\Omega)$ (cf.~\cref{ex:wend}) is generated by Wendland kernels}.
\end{theorem}

We require an auxiliary result based on \cite[Proposition 3.6]{WendRieg05} to prove~\cref{thm:error:regularization}.
\begin{lemma}\label{thm:Wend:regularization}
Let $k\in\N$. Then there are constants $C,h_0 > 0$ such that for every finite set $\calX = \{x_j\}_{j=1}^d \subset \Omega$ of sample points with $h_{\calX} \le h_0$ and all multiindices $\alpha \in \N_0^n$, {$|\alpha|\le k-1$}, we have for all $\lambda \geq 0$, $f\in\calN_{\Phi_{n,k}}(\Omega)$ and $x\in\Omega$
\begin{align*}
\big|D^\alpha f(x) - D^\alpha (R_\calX^\lambda f)(x)\big|\le C\big( h_{\calX}^{k+1/2-|\alpha|} + h_\calX^{-|\alpha|}\sqrt{\lambda}\big)\|f\|_{\calN_{\Phi_{n,k}}(\Omega)}. 
\end{align*}
In particular, for $\alpha=0$,
\begin{align}\label{eq:derivativeconvergence_reg}
|f(x)-R_\calX^\lambda f(x)| \le C\big( h_{\calX}^{k+1/2} + \sqrt{\lambda}\big)\|f\|_{\calN_{\Phi_{n,k}}(\Omega)} \qquad\forall\, f \in \calN_{\Phi_{n,k}}(\Omega),\;x\in\Omega.
\end{align}
\end{lemma}
\begin{proof}
    \cite[Proposition 3.6]{WendRieg05} states that if $\tau = \ell + s$ with $\ell\in\N$, $\ell > n/2$, $s\in (0,1]$, and $0\le j < \ell-n/2$, then $|f - R_\calX^\lambda f|_{W^{j,\infty}(\Omega)}\,\le\,C\Big(h_\calX^{\tau-j-n/2} + h_\calX^{-j}\sqrt\lambda\Big)\|f\|_{H^\tau(\Omega)}$ holds for $f\in H^\tau(\Omega)$, where $|\cdot|_{W^{j,\infty}(\Omega)}$ is the semi-norm $|u|_{W^{j,\infty}(\Omega)} = \sup_{|\alpha|=j}\|D^\alpha u\|_{L^\infty(\Omega)}$.

Let us choose $\tau = \sigma_{n,k} = \frac{n+1}2 + k$. {If $n$ is odd, we have to set $\ell = \tau-1$ and $s=1$. Then $\ell - n/2 = k-1/2 > 0$ so that, if $|\alpha|\le k-1$, setting $j = |\alpha|$ is possible. Hence, $\tau-j-n/2 = k+1/2-|\alpha|$, and the desired estimate follows.} On the other hand, if $n$ is even, then $\ell = k + n/2$ and $s = 1/2$ so that only $0\le j\le k-1$ is allowed. Hence, if $|\alpha|\le k-1$, the result follows.
\end{proof}
\begin{proof}[Proof of~\cref{thm:error:regularization}]
The proof of the theorem follows the lines of that of \cite[Theorem 5.2]{KohnPhil24} together with \cref{thm:Wend:regularization}. %
\end{proof}

{Let us briefly discuss the convergence of $\wh\calK_0$ to $\calK$ in the infinite-data limit, that is, $h_\calX\to 0$ (implying $|\calX|\to\infty$). \cref{thm:error:regularization} establishes 
convergence in the operator norm $\|\cdot\|_{\calN_{\Phi_{n,k}}(\Omega)\to C_b(\Omega)}$ and provides quantitative bounds. However, this result does not imply convergence in operator norms, such as $\|\cdot\|_{C_b(\Omega)\to C_b(\Omega)}$ or $\|\cdot\|_{L^2(\Omega)\to L^2(\Omega)}$, nor does it ensure any form of spectral convergence. For detailed discussions on spectral convergence as well as more sophisticated EDMD variants, see, e.g., \cite{COLBROOK2024127} and~\cite{colbrook2024rigorous}.}

We conclude this section with a remark on alternative choices of kernels.

\begin{remark}\label{rem:Matern}
    The proof of~\cref{thm:error:regularization} on the Koopman approximation consists of two steps. The first ingredient is Koopman invariance of the RKHS, that is, $\mathcal{K}\mathbb{H}\subset \mathbb{H}$. Note that this property has to be checked carefully for chosen kernels. While for Wendland kernels corresponding to fractional Sobolev spaces (and thus also for Matérn kernels generating the same spaces \cite[Example 5.7]{FassYe11}), this invariance holds as shown in \cite{KohnPhil24}, there are also negative results, e.g., for Gaussian kernels~\cite{GonzAbud23}. The second step is establishing an interpolation (or finite-data) error estimate as stated in \cref{thm:Wend:regularization}, e.g., for general radial basis functions in~\cite{narcowich2005sobolev,narcowich2006sobolev}.
\end{remark}

\subsection{Data-driven surrogate dynamics}

In what follows, we propose a Koop\-man-based surrogate model for the dynamics \cref{eq:dynamics}. 
For this, let $\psi_j\in\calN_{\Phi_{n,k}}(\Omega)$, $j\in [1:M]$, be observables such that the map $\Psi = [\psi_1,\ldots,\psi_M]^\top \in\calN_{\Phi_{n,k}}(\Omega)^M$ has a left inverse $\Upsilon : \R^M \to \R^n$, i.e., $\Upsilon(\Psi(x)) = x$ for $x\in\Omega$,
with modulus of continuity $\omega_\Upsilon \in \mathscr{K}$. Note that continuity and injectivity of $\Psi$ imply $M\ge n$. This follows, for example, from the Borsuk-Ulam theorem, see, e.g., \cite[Theorem V.8.9]{Span66}.
{The left inverse map~$\Upsilon$ can be chosen linear if and only if the coordinate projections $x \mapsto x_i$, $i \in [1:n]$, are in the linear span of the lifting~$\Psi$}.

\begin{remark}\label{rem:modulus}
If $\psi_i(x) = x_i$ for $i\in [1:n]$, then $\Upsilon : \R^M\to\R^n$ can be chosen such that $\omega_\Upsilon(r) = r$. Indeed, if $x\in\R^n$ and $y\in\R^{M-n}$, set $\Upsilon(x,y) = x$. Then $\Upsilon(\Psi(x)) = x$ for $x\in\Omega$ and $\|\Upsilon(x,y) - \Upsilon(\bar x,\bar y)\| = \|x-\bar x\|\le \|(x,y) - (\bar x,\bar y)\|$.
\end{remark}

Since the right-hand side of the dynamics \cref{eq:dynamics} satisfies
$$
F(x) = \Upsilon(\Psi(F(x))) = \Upsilon\big(\big[\calK\psi_1(x),\ldots,\calK\psi_M(x)\big]^\top\big),
$$
we may define a data-driven surrogate model by
\begin{align}\label{eq:dynamics:surrogate}
    x^+ = \widehat{F}_\lambda(x) = \Upsilon\big(\big[\widehat\calK_\lambda\psi_1(x),\ldots,\widehat\calK_\lambda\psi_M(x)\big]^{\hspace*{-1mm}\top}\big) = \Upsilon\big(\Psi_{F(\calX)}^\top(\bK_\calX+\lambda I)^{-1}\bk_\calX(x)\big)
\end{align}
with $\Psi_\calZ = [\Psi(z_1),\ldots,\Psi(z_d)]^\top\in\R^{d \times M}$, where we utilize the approximation~$\widehat{K}_\lambda$, $\lambda \geq 0$, of the Koopman operator provided in~\cref{eq:koopman_reg}. 
We note that if the flow map $F$ can be evaluated directly, one may choose $\Psi = \Upsilon = \operatorname{id}$. 

The following result provides an error bound on this approximation by means of the fill distance and the regularization parameter. To this end, set
$$
    \|\Psi\|_{\calN_{\Phi_{n,k}}(\Omega)^M} := \big(\sum\nolimits_{j=1}^M\|\psi_j\|_{\calN_{\Phi_{n,k}}(\Omega)}^2\big)^{1/2}.
$$

\begin{corollary}\label{cor:Fbound_nonproportional}
        Let $k\geq 1$ and $F\in C_b^{\ceil{\sigma_{n,k}}}(\Omega;\R^{n})$. Then there are constants $C,h_0>0$ such that for any finite set $\calX = \{x_i\}_{i=1}^d \subset\Omega$ of sample points with $h_{\calX}\le h_0$ and $\lambda \geq 0$ we have
    \begin{align}\label{eq:Fbound_nonproportional}
        \|F(x) - \widehat{F}_\lambda (x)\| \leq \omega_\Upsilon \Big(C\big(h_\calX^{k+1/2} + \sqrt\lambda\big)\|\Psi\|_{\calN^l_{\Phi_{n,k}}(\Omega)}\Big) \qquad \forall x\in \Omega
    \end{align}
\end{corollary}
\begin{proof}
    Making use of \cref{thm:error:regularization}, the claim follows from the computation
    \begin{align*}
    \|F(x) - \widehat F_\lambda(x)\|
    &= \big\|\Upsilon\big(\big[\calK\psi_1(x),\ldots,\calK\psi_M(x)\big]^\top\big) - \Upsilon\big(\big[\widehat\calK_\lambda\psi_1(x),\ldots,\widehat\calK_\lambda\psi_M(x)\big]^\top\big)\big\|\\
    &\le \omega_\Upsilon\big(\big\|\big[\calK\psi_1(x) - \widehat\calK_\lambda\psi_1(x),\ldots,\calK\psi_M(x) - \widehat\calK_\lambda\psi_M(x)\big]\big\|\big)\\
    &= \omega_\Upsilon\big(\big[\sum\nolimits_{j=1}^M\big|\calK\psi_j(x) - \widehat\calK_\lambda\psi_j(x)\big|^2\big]^{1/2}\big)\\
    &\le \omega_\Upsilon\big(C\left(h_\calX^{k+1/2} + \sqrt{\lambda}\right)\big[\sum\nolimits_{j=1}^M\|\psi_j\|_{\calN_{\Phi_{n,k}}(\Omega)}^2\big]^{1/2}\big).
    \end{align*}
\end{proof}

\section{Data-driven surrogates: Lyapunov stability}
\label{sec:stability}

In this part, we provide our first main result. 
We provide sufficient conditions that ensure transferability of stability results from the original dynamics to the data-driven surrogate and vice versa. 
To be more precise, we show that asymptotic stability of the dynamical system~\cref{eq:dynamics}, certified by a %
Lyapunov function~$V$ with modulus of continuity~$\omega_V$, 
implies semi-global practical asymptotic stability of the kEDMD surrogate model given by \cref{eq:dynamics:surrogate} and vice versa. 
Furthermore, under some compatibility assumptions on the Lyapunov function, semi-global asymptotic stability is even preserved for the unregularized surrogate, i.e., $\widehat{F}_\lambda$ with $\lambda = 0$.
To keep the presentation technically simple(r), in the following we consider the unregularized data-driven surrogate dynamics
\begin{align}\label{eq:dynamics:surrogate:lam0}
   x^+= \widehat F (x) := \widehat F_0(x)
\end{align}
and discuss the regularized surrogate dynamics ($\widehat{F}_\lambda$ with $\lambda > 0$), in subsequent comments. 
Due to the relation $\widehat\calK_0 = \widehat\calK$, this corresponds to vanilla kernel EDMD~\cref{eq:edmd}. 

An important ingredient for the subsequent stability analysis is that equilibria are preserved in the surrogate models if they are contained in the set~$\mathcal{X}$ of data points.
We recall that a state $x^* \in \Omega$ is called an {\em equilibrium} of the dynamics~\cref{eq:dynamics} if $F(x^*) = x^*$. 
\begin{proposition}\label{prop:equilibrium}
    A data point $x^*\in\calX$ is an equilibrium of the dynamics~\cref{eq:dynamics} if and only if it is an equilibrium of the surrogate dynamics~\cref{eq:dynamics:surrogate:lam0}.
\end{proposition}
\begin{proof}
    For each $k\in [1:d]$, we have
    $
        \Psi_{F(\calX)}^\top\bK_\calX^{-1}\bk_\calX(x_k) = \Psi_{F(\calX)}^\top\bK_\calX^{-1}\bK_\calX e_k = \Psi_{F(\calX)}^\top e_k = \Psi(F(x_k)),
    $
    {where $e_k$ denotes the $k$th canonical unit vector}.
    Applying $\Upsilon$, the left inverse of $\Psi$, to both sides of this equation yields $\widehat F(x_k) = F(x_k)$ for all $k\in [1:d]$. Thus, as $x^*\in \calX$, we have $F(x^*) = \widehat F(x^*)$ which implies the claim.
\end{proof}

We point out that the relation of \cref{prop:equilibrium} above is only approximately preserved when considering regularized kEDMD \cref{eq:dynamics:surrogate} with a regularization parameter $\lambda > 0$. More precisely, in view of the implicit function theorem, the data-driven surrogate has an equilibrium in a neighborhood of the original model's equilibrium and the size of this neighborhood is proportional to~${\lambda}$.

We now provide the fundamentals for our subsequent analysis, that is, the notion of a Lyapunov function and (practical) asymptotic stability, as well as a standard result from Lyapunov stability theory.
\begin{definition}\label{def:stability}
    Consider the discrete-time dynamical system~\cref{eq:dynamics} on the set~$\Omega$. 
    Let {$x^*$ be an equilibrium, i.e., $x^* = F(x^*)$, and} let $Y \subset \Omega$ be a forward invariant set w.r.t.\ the dynamics~\cref{eq:dynamics} containing the equilibrium $x^*$ in its interior $\operatorname{int}(Y)$. Then, %
    
    {\rm (i)} the equilibrium~$x^*$ is said to be {\em asymptotically stable} with domain of attraction~$Y$ if there exists $\beta \in \mathscr{K\!L}$ such that 
    \begin{align}\label{eq:stability:asymptotic}
        \| F^n(x) - F^n(x^*) \| \leq \beta(\| x - x^* \|,n)\quad\text{for all $x \in Y$ and $n \in \mathbb{N}_0$.}
    \end{align}
    
    {\rm (ii)} Let, in addition, a set $P \subset Y$ be given that is forward invariant w.r.t.\ the dynamics~\cref{eq:dynamics}. 
    The equilibrium $x^*$ is called $P$-{\em practically asymptotically stable} on $Y$ if there exists $\beta \in \mathscr{K\!L}$ such that, for all $x \in Y$ and $n \in \mathbb{N}_0$, either the inclusion $F^n(x) \in P$ or inequality~\cref{eq:stability:asymptotic} holds.
\end{definition}
{If the set~$P$ is a ball centered at the equilibrium~$x^*$ with radius~$r$, practical asymptotic stability may simply be defined by
\begin{align}
    \| F^n(x) - F^n(x^*) \| \leq \max\{ r, \beta(\| x - x^* \|,n) \} \quad\text{for all $x \in Y$ and $n \in \mathbb{N}_0$.}
\end{align}}%
Asymptotic stability can be characterized by means of Lyapunov functions, see~\cite{HaddChel08}.
\begin{definition}\label{def:lyapunov}
    A continuous function $V : Y \subset \Omega \to \mathbb{R}_{\geq 0}$ is said to be a {\em Lyapunov function} w.r.t.\ the dynamics~\cref{eq:dynamics} and the equilibrium $x^* = F(x^*) \in Y$ if there exist $\alpha_1, \alpha_2 \in \mathscr{K}_\infty$ and $\alpha_V \in \mathscr{K}$ such that
    \begin{align}\label{eq:lyapunov}
        \alpha_1(\| x-x^* \|) \leq V(x) \leq \alpha_2(\|x-x^*\|) \qquad\forall\, x \in Y
    \end{align}
    and the Lyapunov decrease condition given by the inequality
    \begin{align}\label{eq:lyapunov:decrease}
        V(F(x)) \leq V(x) - \alpha_V(\|x-x^*\|) \qquad\forall\,x \in Y \text{ with } F(x) \in Y.
    \end{align}
\end{definition}
The following proposition is assembled from~\cite[Theorem~2.19 and~2.20]{GrunPann17}.
\begin{proposition}\label{prop:lyapunov}
    Let the sets $P,Y \subset \Omega$ satisfying $P \subset Y$ be forward invariant w.r.t.\ the dynamics~\cref{eq:dynamics} with equilibrium $x^* \in \operatorname{int}(P)$. {Then, the following holds:}
    \begin{enumerate}
        \item [(i)] If $V$ is a Lyapunov function in accordance to \cref{def:lyapunov}, the equilibrium~$x^*$ is asymptotically stable on~$Y$.
        \item [(ii)] If $V$ is a Lyapunov function satisfying the decrease condition~\cref{eq:lyapunov:decrease} on $S = Y \backslash P$, $x^*$ is $P$-practically asymptotically stable on $Y$.
    \end{enumerate}
\end{proposition}

\subsection{Inheritance of stability properties}
\label{sec:sub:stability:practical}

Next, using the novel approximation bound of \cref{thm:error:regularization} and the error estimate on the surrogate dynamics~\cref{eq:dynamics:surrogate:lam0}, we present our first main result, i.e., that a Lyapunov function (and, thus, asymptotic stability) of the dynamics~\cref{eq:dynamics} implies practical asymptotic stability w.r.t.\ \cref{eq:dynamics:surrogate:lam0} and vice versa.
Herein, the practical region~$P$ can be rendered to be an arbitrary small neighbourhood of the equilibrium if the fill distance~$h_{\mathcal{X}}$ is sufficiently small. 
\begin{theorem}[Practical asymptotic stability]\label{thm:lyapunov:kEDMD:practical} 
    Let $x^* \in \mathcal{X}$ be an equilibrium w.r.t.\ the dynamics~\cref{eq:dynamics} given by~$F \in \mathcal{C}_b^{\lceil \sigma_{n,k} \rceil}(\Omega; \mathbb{R}^n)$ with {$\sigma_{n,k}$ defined by~\cref{eq:identity:RKHS:Sobolev:Wendland}} and $k \geq 1$ and, thus, also of the data-driven surrogate~\cref{eq:dynamics:surrogate:lam0} represented by~$\widehat{F}$ (see \cref{prop:equilibrium}).
    {Furthermore, let $\psi_i(x) = x_i$, $i \in [1:n]$, be given}. 
    Then, the following two statements hold:
    \begin{itemize}
        \item[{\rm (i)}] {Let $V: \Omega \to \mathbb{R}_{\geq 0}$ be a Lyapunov function w.r.t.\ the dynamics~\cref{eq:dynamics} on~$\Omega$ with modulus of continuity $\omega_V \in \mathscr{K}$. 
            Further, let $c_\Omega > 0$ be sufficiently small such that the sublevel set $V^{-1}(c_\Omega):= \{x \in \Omega \mid V(x) \leq c_\Omega\}$ is closed and the decrease condition~\cref{eq:lyapunov:decrease} holds for all $x \in V^{-1}(c_\Omega)$.\footnote{In particular, $V^{-1}(c_\Omega)$ is, then, forward invariant w.r.t.\ \cref{eq:dynamics}.}}
            Then, $x^*$ is practically asymptotically stable w.r.t.\ \cref{eq:dynamics:surrogate:lam0} in the sense that, for every {$\varepsilon \in (0, \alpha_2^{-1}(c_{\Omega})]$}, the practical region~$P$ can be chosen as a subset of the $\varepsilon$-ball $\mathcal{B}_\varepsilon(x^*)$ if only the fill distance~$h_{\mathcal{X}}$ is sufficiently small.
    \item[{\rm (ii)}] The statement~{\rm (i)} holds upon switching the roles of $F$ and $\widehat F$, {if the comparison functions of the Lyapunov function are uniform for all $h_\mathcal{X}\in (0,h_0]$ for some $h_0$}\footnote{This assumption holds, e.g., for value functions in data-driven predictive control, see \cite{BoldScha25}. 
    The proof follows very similarly to the proof of \cref{lem:Lipschitz:uniform} using the uniform error bound.}, i.e., the existence of a Lyapunov function w.r.t.\ the data-driven surrogate dynamics~\cref{eq:dynamics:surrogate:lam0} {uniform in the fill distance} implies practical asymptotic stability w.r.t.\ the dynamics~\cref{eq:dynamics}.
    \end{itemize}
\end{theorem}
{The imposed smallness condition on the fill distance~$h_\mathcal{X}$ corresponds to a condition on the collected data, which may be easily met by grid-based sampling in the set~$\Omega$ as illustrated in our numerical simulations in~\Cref{subsec:numerics}, but is typically not satisfied when sampling asymptotically stable systems along trajectories.}

{To prove \cref{thm:lyapunov:kEDMD:practical}, we require Lipschitz continuity of the surrogate model, uniformly w.r.t. the approximation accuracy depending on the fill distance $h_\mathcal{X}$ and the regularization parameter $\lambda$. 
To outline a possible extension of \cref{thm:lyapunov:kEDMD:practical} to the approximant~$F_\lambda$ generated with regularization parameter $\lambda > 0$, we show the respective assertion for this more general case. Hereby, it becomes apparent that the regularization parameter~$\lambda$ has to be chosen proportional to the squared fill distance~$h_\mathcal{X}$.
\begin{lemma}[Lipschitz continuity of the data-driven surrogate~$\widehat{F}_\lambda$]\label{lem:Lipschitz:uniform}\ \\
    Let {$k \geq 1$} and $F \in C_b^{\ceil{\sigma_{n,k}}}(\Omega;\R^{n})$ with $\sigma_{n,k} := \frac{n+1}{2} + k$. 
    Then, there exist constants $h_0>0, \bar{L}\geq 0$ such that, 
    if the fill distance~$h_\mathcal{X}$ and the regularization parameter~$\lambda$ satisfy $h_\mathcal{X} \leq h_0$ and $h_\mathcal{X}^{-1} \sqrt{\lambda} \leq c_\lambda$ for some $c_\lambda \in [0,\infty)$, the surrogate~$\widehat{F}_\lambda$ defined by~\eqref{eq:dynamics:surrogate:lam0} is Lipschitz continuous in~$x$ with Lipschitz constant~$L_{\widehat{F}_\lambda} \leq \bar{L}$ , i.e.,
    \begin{align}\label{eq:Lipschitz:surrogate:uniform}
        \| \widehat{F}_\lambda(x) - \widehat{F}_\lambda(y)\| \leq L_{\widehat{F}_\lambda} \|x - y\| \qquad\forall\, x,y \in \Omega.
    \end{align}
\end{lemma}
\begin{proof}    
    $\widehat{F}_\lambda(x) - \widehat{F}_\lambda(y) = [\widehat{F}_\lambda(y) + (x - y)^\top J_{\widehat{F}_\lambda}(\xi)] - \widehat{F}_\lambda(y) = (x - y)^\top J_{\widehat{F}_\lambda}(\xi)$ holds for $x,y \in \Omega$ and some $\xi \in \{(t - 1)x + t y : t \in [0, 1]\}$ by Taylor's theorem, where $J$ stands for the Jacobian.
    Hence, we may estimate
    \begin{align}
        \|\widehat{F}_\lambda(x) \hspace*{-0.5mm}-\hspace*{-0.5mm} \widehat{F}_\lambda(y)\| \leq \|x \hspace*{-0.5mm}-\hspace*{-0.5mm} y\| \|J_{\widehat{F}_\lambda}(\xi) \| \leq \|x \hspace*{-0.5mm}-\hspace*{-0.5mm} y\| \Big[ \|J_{\widehat{F}_\lambda}(\xi) \hspace*{-0.5mm}-\hspace*{-0.5mm} J_{F}(\xi)\| + \|J_{F}(\xi)\| \Big]. \nonumber
    \end{align}
    Since the function~$F$ is continuously differentiable on the compact set~$\Omega$, we have $\| J_{F}(\xi)\| \leq \widetilde{C} \in (0,\infty)$. 
    To compute a bound for $\|J_{\widehat{F}_\lambda}(\xi) - J_{F}(\xi)\|$, we apply \cref{thm:Wend:regularization} with multiindex~$|\alpha| = 1$ to obtain $\| J_{\widehat{F}_\lambda}(z) - J_{F}(z)\| \leq C_{J_{F}} ( h_\mathcal{X}^{k-1/2} + h_\mathcal{X}^{-1} \sqrt{\lambda} )$, 
    where $C_{J_{F}}$ is a constant independent of the fill distance~$h_\mathcal{X}$ and~$\lambda$. 
    Overall, we get the desired inequality~\eqref{eq:Lipschitz:surrogate:uniform} with $L_{\widehat{F}}$ given by $\widetilde{C} + C_{J_{F}} ( h_\mathcal{X}^{k-1/2} + h_{\mathcal{X}}^{-1} \sqrt{\lambda} )$. Hence, the assertion follows with $\bar{L} := \widetilde{C} + C_{J_{F}} ( h_0^{k-1/2} + c_\lambda )$. 
\end{proof}}%
{The continuity in \cref{lem:Lipschitz:uniform} is uniform in the approximation bound~$\varepsilon>0$. This is in contrast to finite-element dictionaries, see~\cite{SchaWort23}, where the derivative of the ansatz functions increases for decreasing mesh size. Here, an error bound on the derivative for kernel-based interpolations is central for proving uniform estimates. Moreover, we note that the adaptive choice of the regularization parameter (depending on the fill distance) stated in \Cref{lem:Lipschitz:uniform} is closely related to the approach of \cite{WendRieg05}.}

{Next, using Lemma~\ref{lem:Lipschitz:uniform}, we prove Theorem~\ref{thm:lyapunov:kEDMD:practical}.}
\begin{proof}
    We begin with assertion~(i) {meaning that $V$ is} a Lyapunov function w.r.t.\ the dynamics~\cref{eq:dynamics}.
    Since $\psi_i(x) = x_i$ for {all} $i\in [1:n]$, it follows that $\Upsilon$ can be chosen such that $\omega_\Upsilon(r)=r$, cf.\ \cref{rem:modulus} {and the subsequent comment}.
    Then, for all $x \in Y:=V^{-1}(c_\Omega)$, we have using \cref{eq:lyapunov:decrease} and \cref{cor:Fbound_nonproportional} constant $C > 0$
    \begin{eqnarray}
        V(\widehat{F}(x)) & = & \Big[ V(\widehat{F}(x)) - V(F(x)) \Big] + V(F(x)) \nonumber \\
        & \leq & \omega_V( \|\widehat{F}(x) - F(x)\| ) + V(x) - \alpha_V(\| x - x^* \|) \nonumber \\
        & \leq & \Big[ \omega_V( \bar{C} h_\mathcal{X}^{k + 1/2} ) - s \alpha_V(\| x - x^* \|) \Big] + V(x) - (1-s)\alpha_V(\| x {- x^*} \|)%
        \label{eq:criticalestimate}
    \end{eqnarray}
    where $s \in (0, 1)$, $(1-s)\alpha_V \in \mathscr{K}$, and $\bar{C} := C \| \Psi \|_{\mathcal{N}^M_{{\Phi_{n,k}}}(\Omega)}$. Hence, if
    \begin{align}\label{eq:lyapunov:proof}
        \omega_V(\bar{C} h_\mathcal{X}^{k+1/2}) \leq s \alpha_V(\| x - x^* \|)    
    \end{align}
    holds, we have a Lyapunov-decrease inequality of the form \cref{eq:lyapunov:decrease} along the data-driven surrogate dynamics~\cref{eq:dynamics:surrogate:lam0}. 
    {Before we proceed, we point out that the inclusion $B_\varepsilon(x^*)\subset V^{-1}(c_\Omega)$ holds due to our choice $\varepsilon \in (0, \alpha_2^{-1}(c_{\Omega})]$: I}f $x \in \mathcal{B}_\varepsilon(x^*)$, then $\|x-x^*\| \leq \alpha_2^{-1}(c_{\Omega})$ and, thus, $\alpha_2(\|x-x^*\|) \leq c_\Omega$. Consequently, $x$ is an element of $V^{-1}(c_\Omega)$ in view of the second inequality in~\cref{eq:lyapunov}. 
    
    We construct a forward invariant set $P\subset B_\varepsilon(x^*)$ such that \cref{eq:lyapunov:proof} holds on $Y\setminus P$. This implies the claim due to \cref{prop:lyapunov}.
     To this end, let $c_\varepsilon > 0$ be such that $P := V^{-1}(c_{\varepsilon}) \subset \mathcal{B}_\varepsilon(x^*)$ holds and choose $\eta > 0$ small enough to ensure the inclusion $\widehat{F}(\mathcal{B}_\eta(x^*)) := \{ \widehat{F}(x) : x\in\mathcal{B}_\eta(x^*)\} \subset P$ using $\widehat{F}(x^*) = x^*$, the {uniform Lipschitz} continuity of the map~$\widehat{F}$ {shown in \cref{lem:Lipschitz:uniform}}, and Condition~\eqref{eq:lyapunov} for the Lyapunov function~$V$. 
    Then, choose $c_\eta > 0$ maximal such that the inclusion $V^{-1}(c_{\eta}) \subset \mathcal{B}_\eta(x^*)$ holds. By definition, for every $x \in V^{-1}(c_{\eta})$, we have $x \in \mathcal{B}_\eta(x^*)$ and, thus, $\widehat{F}(x) \in V^{-1}(c_\varepsilon)$, i.e., $V(\widehat{F}(x)) \leq c_\varepsilon$. Thus, also $\widehat F(x) \in \mathcal{B}_\varepsilon(x^*)$.
    Finally, for every $x \in V^{-1}([c_\eta, c_\Omega])= \{x \in \mathbb{R}^n : c_\eta \leq V(x) \leq c_\Omega\}$, choosing the fill distance~$h_\mathcal{X}$ sufficiently small\footnote{Recall that $\omega_V \in \mathscr{K}$ is strictly increasing. Then, for instance, a fill distance~$h_{\mathcal{X}}$ satisfying $h_\mathcal{X}^{k+1/2} \leq \omega_V^{-1}(s \alpha_V(\|\delta\|)) \bar{C}^{-1} D^{-1}$ with $\delta > 0$ such that $\mathcal{B}_\delta(x^*) \subset V^{-1}(c_\eta)$ and $D$ satisfying $V^{-1}(c_\Omega) \subset \mathcal{B}_D(x^*)$ suffices.}, we get forward invariance of the set~$P$ and the required Lyapunov decrease on $S := Y \setminus P$ showing that $x^*$ is practically asymptotically stable with the practical region~$P$ contained in the desired $\varepsilon$-ball. 

    Next, %
    we show assertion~(ii). If $V$ is a Lyapunov function for the data-driven surrogate dynamics~\cref{eq:dynamics:surrogate:lam0}, the proof follows analogously as in case~(i) with $F$ instead of~$\widehat{F}$. {The uniformity of the comparison functions ensures that the involved quantities do not deteriorate for small fill distances, e.g.,  that \eqref{eq:lyapunov:proof} is satisfied in any neighborhood of the equilibrium when choosing the fill distance small enough.}
\end{proof}
We conclude this subsection by providing two extensions of this result. First, we could have shown semi-global (practical) asymptotic stability: To this end, one needs a global Lyapunov function $V: \mathbb{R}^n \rightarrow \mathbb{R}_{\geq 0}$. Then, for a given compact set~$K$, we determine a sublevel set~$Y$ large enough such that the inclusion $K \subset Y$ holds using that $V$ is proper due to~\cref{eq:lyapunov}. Then, we set $\Omega$ such that $Y \subset \Omega$ and generate the respective data-driven surrogate model. 

Second, the results can be transferred to the dynamics~$\widehat{F}_\lambda$ defined in~\cref{eq:dynamics:surrogate}, resulting from the regularized regression problem ($\lambda > 0$). 
The only change in the proof is the choice of~$\eta$ to ensure $\widehat{F}(\mathcal{B}_\eta(x^*)) \subset P$. 
Here, the regularization parameter has to be chosen sufficiently small such that %
this inclusion is preserved for the given~$\varepsilon$.

\subsection{Proportional bounds and asymptotic stability}
\label{subsec:proportional_interpolation}

We provide conditions under which asymptotic stability is preserved.
First, we present novel proportional bounds on the approximation error extending previous results from~\cite{BoldGrun23,StraScha24:generator}, where comparable, yet probabilistic bounds were obtained under additional assumptions to cover the projection error. %
To this end, we refine~\cref{thm:error:regularization} to reflect the objective of stabilization by means of proportionality of the right-hand side when approaching the equilibrium. 
Then, assuming some compatibility condition linking the modulus of continuity $\omega_V$ and the Lyapunov decrease condition~\cref{eq:lyapunov:decrease}, we show that asymptotic stability is preserved for the data-driven surrogate model and vice versa. 

To this end, a central ingredient will be an inequality controlling the error in the neighborhood of a considered equilibrium $x^*\in\calX$ (i.e., $F(x^*) = x^*)$ of the form
\begin{align}\label{e:proportional_spoiler}
    \|F(x)-\widehat F(x)\| \leq \omega_\Upsilon(\varepsilon \|x-x^*\|)\quad \forall x\in \Omega,
\end{align}
with $\widehat F$ defined in \cref{eq:dynamics:surrogate:lam0}, where $\varepsilon>0$ may be rendered arbitrarily small by using more data points and $\omega_\Upsilon$ denotes the modulus of continuity of the left-inverse $\Upsilon$ of $\Psi$.

In view of the desired bound \cref{e:proportional_spoiler}, observe that the previously established estimate of \cref{cor:Fbound_nonproportional} is sub-optimal in the case $\la=0$: For data points $x\in\calX$, the left-hand side in the estimate \cref{eq:Fbound_nonproportional} vanishes while the right-hand side remains constant. Thus, the next theorem improves the error bounds of \cref{cor:Fbound_nonproportional} by including a dependence on the distance to the data set in the upper bound. In its proof and in the sequel, we work with the spaces of continuous and bounded functions with continuous and bounded derivatives defined as follows. 
For $s \in \N$, we consider the space $C_b^s(\Omega)$ of continuous functions $f : \Omega\to\R$ for which $D^\alpha f$ is bounded on $\Omega$ for all multiindices~$\alpha \in \N_0^d$, $|\alpha| \leq s$, endowed with the norm
   $ \|f\|_{C_b^s(\Omega)} := \sum\nolimits_{|\alpha|\leq s} \,\sup\nolimits_{x\in \Omega} |D^\alpha f(x)| = \sum\nolimits_{|\alpha|\leq s} \|D^\alpha f\|_{C_b(\Omega)}.
$
\begin{theorem}[Proportional error bounds]\label{thm:bounds:proportional}
    Assume that $k \geq 1+\frac{(-1)^n+1}{2}$. Then, there exist constants $C, h_0 > 0$ such that, for every set $\calX = \{x_j\}_{j=1}^d \subset \Omega$ of sample points satisfying the condition $h_{\calX} \leq h_0$, we have
    \begin{align*}
        \big| f(x) - (P_\calX f)(x) \big| \leq C h_{\calX}^{k-1/2} \dist(x,\calX) \|f\|_{\calN_{\Phi_{n,k}}(\Omega)} \qquad\forall\, f \in \calN_{\Phi_{n,k}}(\Omega), x\in\Omega.
    \end{align*}
    Further, if additionally $F \in C_b^{\ceil{\sigma_{n,k}}}(\Omega;\R^{n})$ for~\cref{eq:dynamics}, there are constants $C, h_0>0$ such that, for any finite set $\calX = \{x_i\}_{i=1}^d \subset \Omega$ of sample points with $h_{\calX}\le h_0$, we have
    \begin{align*}    
        |(\mathcal{K} f)(x) - (\widehat{\mathcal{K}} f)(x)| \leq Ch_\calX^{k-1/2}\dist(x,\calX)\|f\|_{\calN_{\Phi_{n,k}}(\Omega)} \qquad\forall f \in \calN_{\Phi_{n,k}}(\Omega), x\in\Omega.
    \end{align*}
\end{theorem}
\begin{proof}
    First, w prove the first inequality w.r.t.\ interpolation errors. %
    For given $f \in \calN_{\Phi_{n,k}}(\Omega)$, we define the error function by $e_f(x) := f(x) - (P_\calX f)(x)$, $x\in\Omega$. 
    Since $k \ge 1$, it follows from Sobolev's embedding theorem (see, e.g., \cite[Theorem 4.12]{Adams2003}) that $\calN_{\Phi_{n,k}}(\Omega) = H^{\sigma_{n,k}}(\Omega)\subset C_b^1(\Omega)$ and hence $e_f \in C_b^1(\Omega)$. 
    {Using $e_f(x_j) = 0$}, we infer $|e_f(x)| = |e_f(x) - e_f(x_j)| \le (\sup_{z\in\Omega}\,\|\nabla e_f(z)\|) \| x-x_j \|$ for arbitrary $j \in [1:d]$.
    As this holds for all $j \in [1:d]$, we conclude $|e_f(x)| \le (\sup_{z\in\Omega}\,\|\nabla e_f(z)\|) \dist(x,\calX)$. 
    Now, from \cref{thm:Wend:regularization} with $\lambda = 0$ and $|\alpha| = 1$ corresponding to all first derivatives, we obtain $\|\nabla e_f(z)\|^2 = \sum_{i=1}^{n}|\partial_ie_f(z)|^2\,\le\,nC^2h_\calX^{2k-1}\|f\|_{\calN_{\Phi_{n,k}}(\Omega)}^2$ for all $z \in\Omega$. The claim follows after redefining $C:=\sqrt{n}C$.

    Next, we show the second claim.
    Due to \cite[Theorem 4.2]{KohnPhil24}, $F\in C_b^{\ceil{\sigma_{n,k}}}(\Omega;\R^{n})$ implies the invariance of the native space $\calN_{\Phi_{n,k}}(\Omega)$ under the Koopman operator $\mathcal{K}$, and the restricted operator $\mathcal{K} : \calN_{\Phi_{n,k}}(\Omega)\to\calN_{\Phi_{n,k}}(\Omega)$ is bounded. Therefore, for $f\in\calN_{\Phi_{n,k}}(\Omega)$ and $x\in\Omega$, we use the just established result to estimate to prove the claim:
    \begin{align*}
        | (\mathcal{K} f)(x) - (\widehat{\mathcal{K}} f)(x)|
    &= |(\mathcal{K} f)(x) - (P_\calX\mathcal{K} f)(x)|\,\le\, Ch_\calX^{k-1/2}\dist(x,\calX)\|\mathcal{K} f\|_{\calN_{\Phi_{n,k}}}\\
    &\le C\|\mathcal{K}\|_{\calN_{\Phi_{n,k}}\to\calN_{\Phi_{n,k}}}h_\calX^{k-1/2}\dist(x,\calX)\|f\|_{\calN_{\Phi_{n,k}}}.
    \end{align*}
\end{proof}

We now provide a proportional error bound for the surrogate dynamics \cref{eq:dynamics:surrogate:lam0} as claimed in \cref{e:proportional_spoiler}. 

\begin{corollary}\label{c:surrogate_proportional}
    Under the assumptions of \cref{thm:bounds:proportional}, there are constants $C,h_0>0$ such that for any finite set $\calX = \{x_i\}_{i=1}^d \subset\Omega$ of sample points with $h_{\calX}\le h_0$ we have $\|F(x) - \widehat F(x)\|\le \omega_\Upsilon\big(Ch_\calX^{k-1/2}\|\Psi\|_{\calN_{\Phi_{n,k}}(\Omega)^M}\dist(x,\calX)\big)$ and, in particular, 
    \begin{align*}    
        \|F(x) - \widehat F(x)\|\le \omega_\Upsilon\big(Ch_\calX^{k-1/2}\|\Psi\|_{\calN_{\Phi_{n,k}}(\Omega)^M}\|x-x^*\|\big) \qquad\forall\, x^* \in \calX.
    \end{align*}
    \end{corollary}
\begin{proof}
    The first claim follows from the proof of \cref{cor:Fbound_nonproportional} using \cref{thm:bounds:proportional} in the second-last estimate, while the second directly from the definition of the distance.
\end{proof}

Next, we leverage the derived proportional bounds on the approximation error in order to show that asymptotic stability is preserved provided that either a compatibility condition is satisfied or the Lyapunov function is given some norm.
\begin{theorem}[Asymptotic stability]\label{thm:lyapunov:kEDMD:stability}
    Let $x^* \in \mathcal{X}$ be an equilibrium w.r.t.\ the dynamics~\cref{eq:dynamics} given by~$F \in \mathcal{C}_b^{\lceil \sigma_{n,k} \rceil}(\Omega; \mathbb{R}^n)$ with {$\sigma_{n,k}$ defined by~\cref{eq:identity:RKHS:Sobolev:Wendland}} and $k \geq 1$ and, thus, also of the data-driven surrogate~\cref{eq:dynamics:surrogate:lam0} represented by~$\widehat{F}$.\footnote{See \cref{prop:equilibrium}.}
    {Furthermore, let $\psi_i(x) = x_i$, $i \in [1:n]$, be given. 
    In addition, suppose} that one of the following conditions on the function $V: \Omega \to \mathbb{R}_{\geq 0}$ hold:
    \begin{itemize}
        \item[{\rm (a)}] $V$ has a modulus of continuity $\omega_V \in \mathscr{K}$ satisfying the compatibility condition\footnote{{Tacityly, 
            we assume $\omega_V(r) = r^a + \mathcal{O}(r^b)$ and $\alpha_V(r) = r^{\tilde{a}} + \mathcal{O}(r^{\tilde{b}})$ with $b > a > 0$ and $\tilde{b} > \tilde{a} > 0$. 
            Then, condition~\cref{eq:compatability} corresponds to $a \geq \tilde{a}$ and is equivalent to the more general condition $\exists\veps>0 : \limsup_{r\searrow 0}\ \frac {\omega_V(\veps r)}{\alpha_V(r)} < 1$, which will be used in proof and has to be directly verified if the tacitly assumed representation does not hold}.}
        \begin{align}\label{eq:compatability}
            \limsup_{r\searrow 0}\ \frac {\omega_V(r)}{\alpha_V(r)} < \infty.
        \end{align}
    \item[{\rm (b)}] $V(x) = \| x - x^* \|^p$ for some $p\in\N$, and {$\alpha_V(r) = c_V r^p$}.\footnote{
        {Clearly, this implies also the compatibility condition~\cref{eq:compatability}}.} 
    \end{itemize}
    Further, we define the sublevel set $V^{-1}(c_\Omega):= \{x \in \mathbb{R}^n \ | \ V(x) \leq c_\Omega\}$, where $c_\Omega > 0$ is chosen such that $V^{-1}(c_\Omega) \subset \Omega$ is closed. 
    Then, the following two statements hold.
    \begin{itemize}
        \item[{\rm (i)}] Let $V$ be a Lyapunov function w.r.t.\ the dynamics~\cref{eq:dynamics} on~$\Omega$ such that the decrease condition~\cref{eq:lyapunov:decrease} holds for all $x \in V^{-1}(c_\Omega)$, which renders $V^{-1}(c_\Omega)$ forward-invariant. %
        Then, the equilibrium~$x^*$ is asymptotically stable w.r.t.\ the data-driven surrogate~\cref{eq:dynamics:surrogate} for sufficiently-small fill distance~$h_{\mathcal{X}}$.
        \item[{\rm (ii)}] The statement of assertion~(i) holds upon switching the roles of $F$ and $\widehat F$, {if the comparison functions of the Lyapunov function are uniform for all $h_\mathcal{X}\in (0,h_0]$ for some $h_0$}, i.e., the existence of a Lyapunov function w.r.t.\ the data-driven surrogate dynamics~\cref{eq:dynamics:surrogate} {uniform in the fill distance} implies asymptotic stability w.r.t.\ the dynamics~\cref{eq:dynamics} meaning that asymptotic stability is preserved.
    \end{itemize}
\end{theorem}
\begin{proof}
    We begin with assertion~(i). Let $V$ be a Lyapunov function w.r.t.\ \cref{eq:dynamics} satisfying {condition}~(a).
    We proceed along the lines of the proof of \cref{thm:lyapunov:kEDMD:practical}.
    Using the compatibility condition~\cref{eq:compatability} and the proportional bound of \cref{c:surrogate_proportional} in the estimate \eqref{eq:criticalestimate}, we arrive at a counterpart of inequality~\cref{eq:lyapunov:proof} 
    with proportional left-hand side, i.e., $\omega_V(\bar{C} h_\mathcal{X}^{k-1/2}\|x-x^*\|) \leq s \alpha_V(\| x - x^* \|)$ holds for all $x\in V^{-1}(c_\Omega)$ for sufficiently small fill distance. This allows to infer asymptotic stability of the equilibrium~$x^*$ w.r.t.\ the data-driven surrogate dynamics~\cref{eq:dynamics:surrogate} using \cref{prop:lyapunov}. 

    Let us now assume that (b) holds. 
    Then, {using $J_F$ to denote the Jacobian of~$F$},
    \begin{eqnarray*}
        V(\widehat{F}(x)) & = & %
        \| \widehat{F}(x) - F(x) + F(x) - x^* \|^p \leq \Big( \| \widehat{F}(x) - F(x) \| + \| F(x) - x^* \| \Big)^p \\
        & = & \sum\nolimits_{j=0}^p \binom{p}{j} \| \widehat{F}(x) - F(x) \|^{p-j} \| F(x) - x^* \|^j.\\
        & \leq & \|F(x) - x^*\|^p + \sum_{j=0}^{p-1} \binom{p}{j} \left( C h_\calX^{k-1/2}\|x-x^*\| \right)^{p-j} \max_{\xi\in\Omega}\| J_F(\xi)\|^j \|x-x^*\|^j\\
        & \leq & {V(F(x)) + h_\calX^{k- 1/2} \hspace*{-1mm} \underbrace{\left[ \sum_{j=0}^{p-1} \binom{p}{j} C^{p-j} h_0^{(k- \frac 12)(p-j-1)}\max_{\xi\in\Omega}\| J_F(\xi) \|^j \right]}_{=: \bar{C}} \hspace*{-1mm} \|x-x^*\|^{p}} \\
        & \stackrel{\cref{eq:lyapunov:decrease}}{\leq} & V(x) - \alpha_V(\|x-x^*\|) + h_\calX^{k-1/2} {( {\bar{C}}/{c_V}) \alpha_V(\|x-x^*\|)} \\
        & \leq & V(x) - { \Big(1 - h_\calX^{k-1/2} \bar{C} / c_V) \Big) \cdot \alpha_V(\|x-x^*\|)},
    \end{eqnarray*}
    where {$s\in (0,1)$.
    Thus, the assertion follow}s for sufficiently small fill distances~$h_{\mathcal{X}}$. 
    
    {The second statement (ii) follows, analously as in \cref{thm:lyapunov:kEDMD:practical}, upon switching the roles of $\hat F$ and $F$ using uniformity of the Lyapunov function for small fill distances.}
\end{proof}

\subsection{Example and numerical simulations}
\label{subsec:numerics}

We illustrate our findings by an instructive example. Consider the nonlinear discrete-time system from~\cite{liHarf14} given by
\begin{align}\label{eq:Kellett}
    x^+ = F(x) := \frac 18
    \begin{pmatrix}
        \|x\|^2-1 & -1 \\
        1 & \|x\|^2-1 
    \end{pmatrix}x
\end{align}    
with state~$x = (x_1, x_2)^\top$ in the compact set $\Omega = [-2,2]^2$, for which $V(x) = \|x\|^2$ serves as quadratic Lyapunov function on~$\Omega$ w.r.t.\ the equilibrium $x^* = 0$, see~\cite{liHarf14}. 
We use the coordinate maps $\psi_i(x) = x_i$, $i \in \{1, 2\}$, as observables in~\cref{eq:dynamics:surrogate} and set the smoothness degree $k = 1$ of the Wendland kernels. As the Lyapunov function is given by the Euclidean norm, we can directly apply \cref{thm:lyapunov:kEDMD:stability} using assumption (b) to conclude asymptotic stability of the data-driven surrogate dynamics~\cref{eq:dynamics:surrogate:lam0} and vice versa. %
When approximating $F$ by means of a regularized surrogate ($\lambda > 0$), we apply \cref{thm:lyapunov:kEDMD:practical} (and the comments thereafter) to deduce practical asymptotic stability.

\begin{figure}[htb]
    \centering
        \includegraphics[width=0.9 \linewidth]{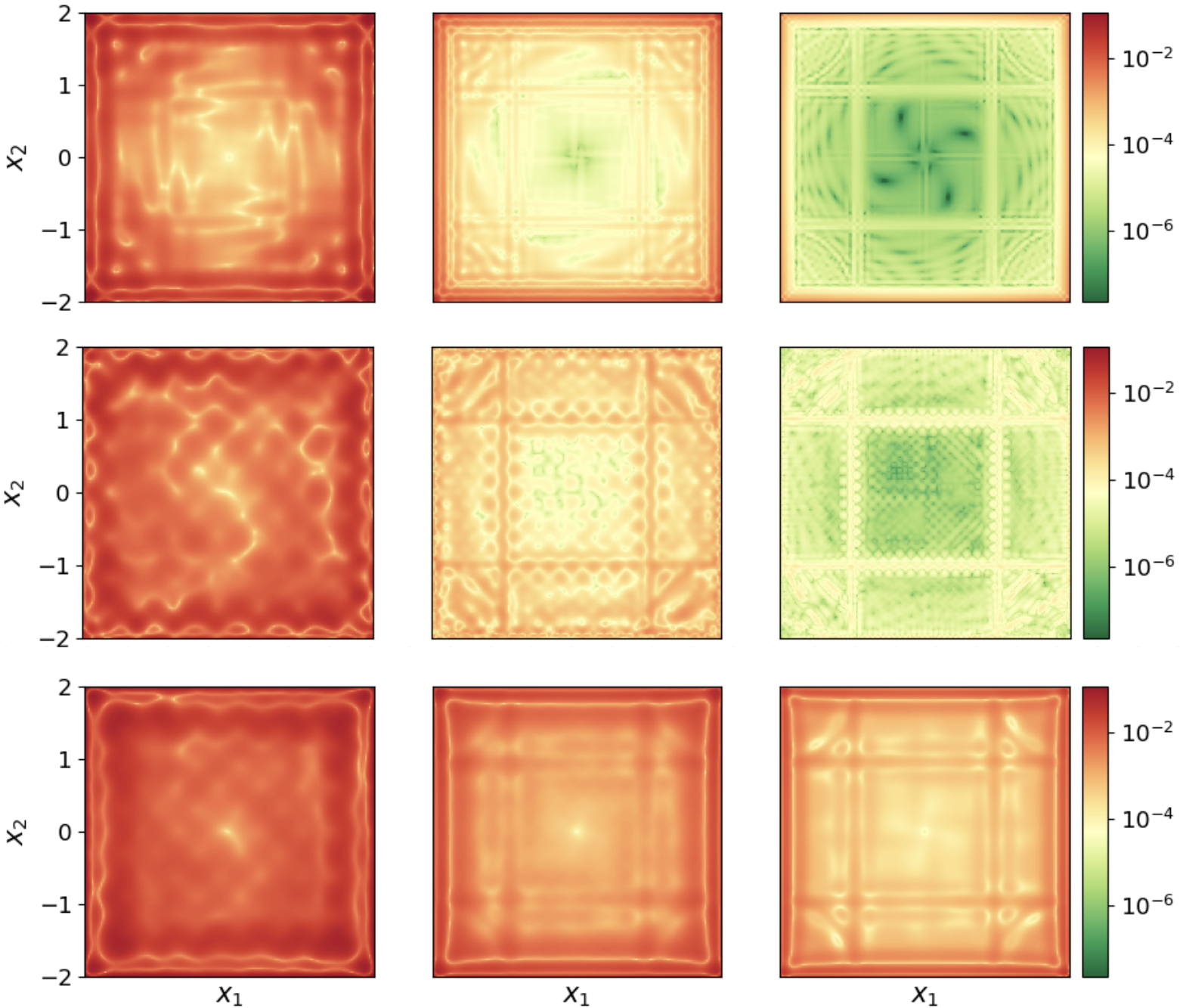}\\[-2mm]
    \caption{Approximation error $\|\widehat{F}(x) - F(x)\|_2$ for system~\cref{eq:Kellett} with $d \in \{441, 1681, 6561\}$ data points for a uniform \braces{mesh sizes~$\delta \in \{0.2, 0.1, 0.05\}$; top} and { a Padua grid \braces{middle}, $d \in \{436, 1654, 6556\}$}. Bottom: Error $\|\widehat{F}_\lambda(x) - F(x)\|_2$ with regularization parameter~$\lambda = 0.01$.}
    \label{fig:heatplot}
\end{figure}

To validate our findings, we conduct numerical computations using a uniform grid $\calX = \delta \mathbb{Z}^2 \cap \Omega$ with mesh size~$\delta > 0$ leading to a fill distance $h_{\calX} = 2^{-1/2}\delta$. 
In addition, we consider a grid consisting of %
{Padua points as described} in~\cite{Cali05}. {For degree~$\gamma_{\text{P}}\in \N$, the Padua points of family~$s = 1$ are given by the set of points
\begin{align*}
    \text{Pad}_{\gamma_P} = \Big\{ P \left(\frac{k\pi}{\gamma_{\text{P}}(\gamma_{\text{P}} + 1)}\right) \mid k = 0, \dots, \gamma_{\text{P}}(\gamma_{\text{P}} + 1)\Big\},
\end{align*}
where~$P:[0, \pi] \rightarrow \R^2$ is defined by 
    $P(t) = \left( -\cos((\gamma_{\text{P}} + 1)t) \ -\cos(\gamma_{\text{P}}t) \right)^\top$.
This results in $|\text{Pad}_{\gamma_P}| = (\gamma_P + 1)(\gamma_P + 2)/ 2$ grid points. Using a grid of Padua points reduces the condition number for two-dimensional problems by minimizing the Lebesgue constant, similar to Chebychev nodes in one dimension.} %
First, in \cref{fig:heatplot}, we inspect the one-step prediction error %
of the surrogate model~$\widehat{F}$ without regularization term, i.e.\ $\lambda = 0$. 
The (uniform) validation grid $\calX_{\text{val}}$ uses the mesh size~$\delta_{\text{val}} = 0.025$ and is given by $\calX_{\text{val}} = \left(\delta_{\text{val}}\mathbb{Z} + \frac{\delta_{\text{val}}}{2}\right)^2 \cap \Omega$. 
In the top row of \cref{fig:heatplot}, the results for the uniform grid with  mesh size~$\delta\in \{0.2 \cdot 2^{-i}\,|\,i \in [0:2]\}$ (and, therefore, $h_\calX \in \{0.2 \cdot 2^{-i-\frac{1}{2}}\,|\,i \in [0:2]\}$) can be seen. 
The intensity plots of the error~$\|\widehat{F}(x) - F(x)\|_2$ for each point~$x$ in the validation grid~$\calX_{\text{val}}$ %
show that the error decreases for decreasing fill distance.  %
The lower row of \cref{fig:heatplot} depicts the results for the associated {Padua-based mesh with degrees~$\gamma_{\text{P}} \in \{28, 56, 113\}$, which results in $d \in\{435, 1653, 6555\}$ grid points.}
Again, we observe that the error decreases the more data points are chosen (as, correspondingly, the fill distance is decreased). In particular, the choice of a {Padua-based} %
grid alleviates large errors at the boundary. 
This also effect the maximum errors $\max_{x \in \calX_{\text{val}}} \|\widehat{F}(x) - F(x)\|_2$, see \cref{tab:max}. %
{On the full domain, the usage of the Padua grid reduces the maximal error visibly.} The closer we get to the origin of the domain, the smaller the maximal error, which indicates the proportional decrease proven in \cref{cor:Fbound_nonproportional}. It can also be seen that for smaller validation areas, i.e.\ the further afar the validation area is from the boundary, the better the surrogate model~$\widehat{F}$ on the uniform grid performs, {and produces errors comparable to the ones obtained from the Padua-based} grid. 
\begin{table}[htb]
    \centering
    \begin{tabular}{|c|c|c|c|c|c|c|}
    \hline
        & {\footnotesize $d=441$} & {\footnotesize $d=435$} & {\footnotesize $d=1681$} & {\footnotesize $d=1653$} & {\footnotesize $d=6561$} & {\footnotesize $d=6555$} \\
       & {\footnotesize uniform} & {\footnotesize Padua} & {\footnotesize uniform} & {\footnotesize Padua} & {\footnotesize uniform} & {\footnotesize Padua} \\
       \hline 
        $[-2.0, 2.0]^2$ & 0.1205 & {0.\textbf{0}127} & 0.\textbf{0}3770 & {0.\textbf{000}79} & 0.\textbf{00}9540 & {0.\textbf{000}150}\phantom{0} \\
        $[-1.0, 1.0]^2$ & 0.\textbf{00}53 & {0.\textbf{00}44} & 0.\textbf{000}30 & {0.\textbf{000}33} & 0.\textbf{0000}21 & {0.\textbf{0000}38}\phantom{0} \\
        $[-0.5, 0.5]^2$ & 0.\textbf{000}7 & 0.{\textbf{000}8} & 0.\textbf{0000}4 & {0.\textbf{0000}2}%
        & 0.\textbf{00000}1 & {0.\textbf{000000}9}%
        \\ %
         \hline
    \end{tabular} 
    \caption{Maximal error for one-step prediction of System~\cref{eq:Kellett} using the kEDMD surrogate~$\widehat{F}$ for uniform and {Padua-based}
    grid on different snippets of the validation grid.} 
    \label{tab:max}
\end{table}
Next, we inspect regularized kEDMD with regularization parameter~$\lambda  = 0.01$. %
{While the results resemble the ones %
depicted in \cref{fig:heatplot} without regularization, the error is increased by about two orders of magnitude close to the origin %
due to the additional term $\sqrt{\lambda}$ in %
\cref{cor:Fbound_nonproportional}}. 

\begin{figure}[htb]
    \centering
    \includegraphics[height=0.3\linewidth]{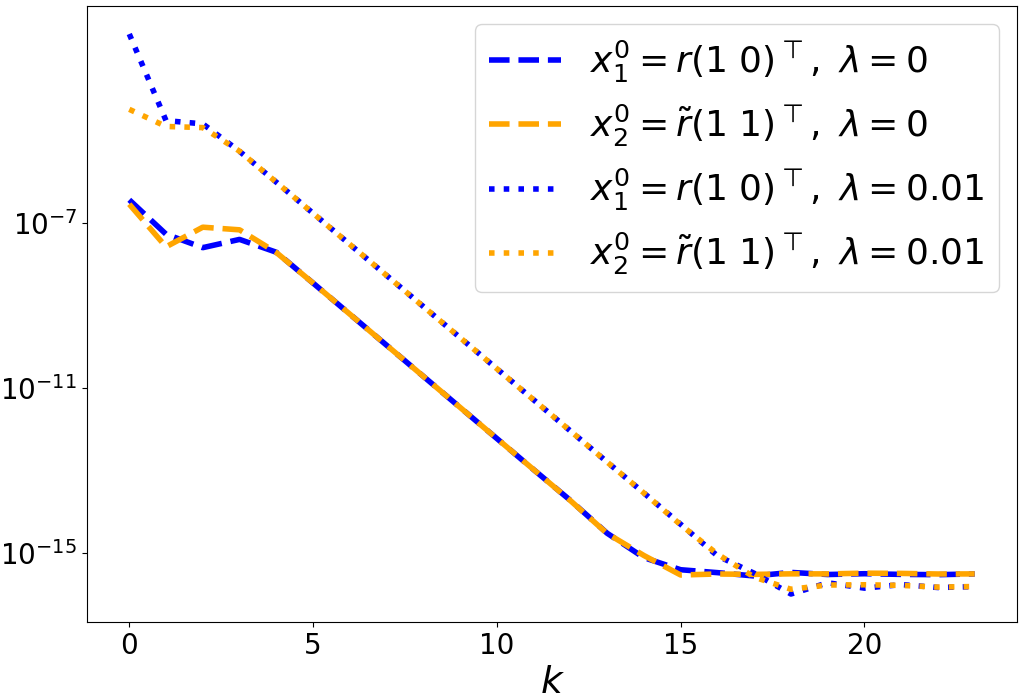}\hspace*{1cm}
    \includegraphics[height=0.3\linewidth]{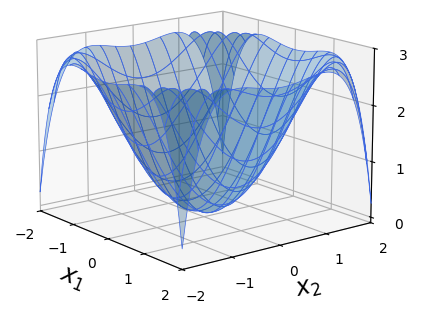}
    \caption{Left: %
    $\| \widehat{F}_\lambda(x(k; x^0_j)) - F(x(k; x^0_j)) \|_2$ along trajectories of~\cref{eq:Kellett} for $r = 1.9$, $\tilde{r} = 0.95 \sqrt{2}$. %
    Right: Lyapunov decrease $V(x)-V(\widehat{F}(x)) - \alpha_V(\|x\|)$ for system~\cref{eq:Kellett} with the kEDMD surrogate~$\widehat{F}$.}
    \label{fig:multi-step}
\end{figure}

For a long-term evaluation of the errors, we choose a uniform grid~$\calX = \delta \mathbb{Z}^2 \cap \Omega$ with $\delta = 0.2 \cdot 2^{-3} = 0.025$ to learn the surrogate dynamics~$\widehat{F}_\lambda$ for the regularization parameters~$\lambda = 0$ and $\lambda = 0.01$.
For the two initial conditions $x^0_1 = r (1, 0)^\top$ and $x^0_2 = r (\cos(\pi/4 ), \sin(\pi/4) )^\top$ that are located on the circle centered at the origin with radius $r = 1.9$ the one-step errors~$\| \widehat{F}_\lambda(x(k; x^0_j)) - F(x(k; x^0_j)) \|_2$ for $k \in [0:20]$ and $j \in \{1, 2\}$ are visualized in \cref{fig:multi-step} in the image on the left.
The difference of the one-step errors decays along the asymptotically stable trajectory of the original dynamics, confirming the proportional error bound \cref{c:surrogate_proportional} and the asymptotic stability of the surrogate model \cref{thm:lyapunov:kEDMD:stability}. 
Lastly, we validate the Lyapunov decrease condition~\cref{eq:lyapunov:decrease}, i.e., $V(\widehat{F}(x)) - V(x) + \alpha_V(\|x\|) \leq 0$, where the data-driven surrogate dynamics~$\widehat{F}$ of \cref{eq:Kellett} are generated using the mesh size~$\delta = 0.05$. In the right plot of \cref{fig:multi-step}, the value of $-(V(\widehat{F}(x)) - V(x) + \alpha_V(\|x\|))$ for $x \in \calX_{\text{val}}$ is plotted and we can observe that the decrease condition is in fact preserved as stated in~\cref{thm:lyapunov:kEDMD:stability}.  

\section{Koopman approximants for control systems}\label{sec:kEDMD:control}

In this section, we consider the discrete-time control-affine system 
\begin{align}\label{eq:dynamics:control}
    x^+ = f(x,u) = g_0(x) + G(x)u,
\end{align}
where $x\in\Omega\subset\R^n$ and the control $u$ is contained in a bounded set $\bbU\subset\R^m$, $m\le n$, with $\linspan\bbU = \R^m$. Here, we assume that the maps $g_0 :\Omega\to\R^n$ and $G : \Omega\to\R^{n\times m}$ are locally Lipschitz continuous, i.e., there is $L_g > 0$ and $L_G > 0$ such that
$$
    \|g_0(x)-g_0(y)\|\le L_g\|x-y\| \qquad\text{and}\qquad
    \|G(x)-G(y)\|\le L_G\|x-y\| \qquad\forall\,x,y\in\Omega
$$
Further, all entries of both $g_0$ and $G$ are assumed to be in the RKHS $\bbH = \calN_{\Phi_{n,k}}(\Omega)$ for some fixed $k\in\N$. Note that, in view of \cref{eq:identity:RKHS:Sobolev:Wendland}, $\calN_{\Phi_{n,k_1}}(\Omega) \subset\calN_{\Phi_{n,k_2}}(\Omega)$ for $k_1\geq k_2$.
\begin{remark}[Continuous-time control-affine dynamics]
    Often, systems~\cref{eq:dynamics:control} are derived from continuous-time control-affine systems governed by the dynamics
    \begin{align}\label{eq:ode:control}
        \dot{x}(t) := \tilde{f}(x(t),u(t)) := \tilde{g}_0(x(t)) + \sum\nolimits_{i=1}^{m} u_i(t)\tilde{g}_i(x(t)) 
    \end{align}
    with locally Lipschitz maps $\tilde{g}_i: \Omega \rightarrow \mathbb{R}^n$, $i \in [0:m]$. 
    Similarly to \cref{rem:ode}, we obtain $x(\Delta t;\bar{x},u) = \bar{x} + \int\nolimits_0^{\Delta t} \tilde{g}_0(x(s;\bar{x},u))\,\mathrm{d}s + \int_0^{\Delta t} \tilde{G}(x(s;\bar{x},u))u(s)\,\mathrm{d}s$
    for a time step $\Delta t > 0$, 
    where $\tilde{G}$ is defined analogously to~$G$ above --~again tacitly assuming existence and uniqueness of the solution $x(\cdot;\bar{x},u)$ on $[0,\Delta t]$. 
    Then, for constant 
    control function~$u$ %
    on the interval $(0,\Delta t)$, a Taylor-series expansion of the solution at $t=0$ yields $x(\Delta t;\bar{x},u) = \bar{x} + \Delta t [ \tilde{g}_0(\bar{x}) + \tilde{G}(\bar{x})u ] + \mathcal{O}(\Delta t^2)$, where we have invoked compactness of~$\mathbb{U}$ and Lipschitz continuity of $f(\cdot,u)$ and $\frac \partial {\partial x}f(\cdot,u)$ on $\overline{\Omega}$ identifying the (constant) control functions with the respective control value. %
    Thus, we obtain a discrete-time system~\cref{eq:dynamics:control} up to an arbitrarily small error if the time step~$\Delta t$ is sufficiently small.
\end{remark}

In this section, we propose a novel learning architecture to learn Koopman approximants of the control-affine system~\cref{eq:dynamics:control}, that allows for flexible sampling of state-control pairs $(\bar{x}_i, \bar{u}_i, \bar{x}_i^+)$ with $\bar{x}_i \in \Omega$, $\bar{u}_i \in \bbU$, and $\bar{x}_i^+ = f(\bar{x}_i,\bar{u}_i)$, $i\in [1:d]$.
The result of our proposed kEDMD-based algorithm are control-affine surrogate dynamics
\begin{align}\label{eq:dynamics:control_surr}
    x^+ = \widehat f(x,u) = \widehat g_0(x) + \widehat G(x)u
\end{align}
using regression and interpolation in Wendland native spaces. In \cref{t:control_main}, we prove a uniform bound on the approximation error $\|f(x,u) - \wh f(x,u)\|$, similarly as in \cref{cor:Fbound_nonproportional} for the approximants $\widehat{F}$ and $\widehat F_\lambda$ of %
system~\cref{eq:dynamics}.
Then, we exemplarily apply the derived proportional bound for feedback stabilization.
To be more precise, we prove in \Cref{subsec:feedback} that if a feedback stabilizes the surrogate dynamics~\cref{eq:dynamics:control_surr}, the desired set point is stabilized in the original system (and vice versa). %

\subsection{Koopman approximants, flexible sampling, and error bounds}

We propose a novel learning architecture to generate Koopman approximants for control-affine systems with the following two key features: On the one hand, the algorithm can be used with 
(almost) arbitrary state-control data allowing for highly flexible sampling including i.i.d., grid-based and trajectory-based data generation. 
On the other hand, we rigorously show bounds on the full approximation error with explicit convergence rates in the infinite-data limit\footnote{{At this point, the choice of the sampling scheme becomes crucial: while small fill distances can be readily achieved using grid-based methods, trajectory-based sampling typically relies on ergodicity---a restrictive property when considering control applications.}} and, foremost, allow for controller design with closed-loop guarantees without imposing restrictive invariance assumptions on the dictionary in EDMD.
To the best of our knowledge, this combination is, up to now, unique. Koopman-based approximations of control-affine systems in an RKHS framework with flexible sampling were also presented in~\cite{bevanda2024nonparametric}, where the kernel (and thus also the RKHS) is suitably extended to capture the control dependency. Therein, however, no error bounds were provided. In addition, our novel approach allows to counteract numerical ill-conditioning by decomposing the approximation process in two steps labeled as macro and micro level, see the following paragraph for details.

The method proposed in \cref{alg:cap} consists of two steps. 
First, we determine $d$~clusters of size~$N \geq m+1$. %
Each cluster corresponds to a center $x_\ell \in \Omega$, $\ell \in [1:d]$, and its $N$ nearest neighbors. 
The chosen centers~$x_\ell$, $\ell \in [1:d]$, determine the fill distance~$h_\mathcal{X}$. 
Then, we use the data triples~$(\bar{x}_i,\bar{u}_i,\bar{x}_i^+)$ corresponding to the cluster centered at~$x_\ell$ to approximate the function values $g_0(x_\ell)$ and $G(x_\ell)$ for each $\ell \in [1:d]$. {This step addresses the absence of direct measurements of $g_0$ and $G$, as these quantities are only available through the lens of the dynamics via a linear combination scaled with the control values.}
In a second step, %
we apply RKHS-based interpolation to approximate the entries in the vector- and matrix-valued functions $g_0$ and $G$ in the $\sup$-norm. We briefly highlight that the clustering step in particular alleviates the inherent ill-conditioning of kernel-based interpolation tasks, which manifests in the rapidly decaying eigenvalues of the kernel matrix~$\bK_\mathcal{X}$ \cite[Theorem 4.16]{Sant18}.
This feature may be seen in the estimate of the subsequent~\cref{t:control_main} as the term including the inverse of the kernel matrix may be controlled by the cluster radius $\varepsilon$.

\begin{algorithm}[htb]
\caption{Data-based approximation of the controlled flow map}\label{alg:cap}

\smallskip
{\it Input:} Triples $(\bar x_i,\bar u_i,\bar x_i^+)$, $i\in [1:\bar d]$, where $\bar x_i^+ = f(\bar x_i,\bar u_i)$, $\bar x_i\in\Omega$, $\bar u_i\in\bbU$, $i\in [1:\bar d]$. Number of clusters~$d \geq \lfloor \bar{d}/N \rfloor$ and minimal number of cluster elements $N \geq m + 1$.

\smallskip\noindent
{\it Output:} Approximation $\wh f$ as in \cref{eq:dynamics:control_surr}. %

\smallskip\hrule
\smallskip

\textit{Initialisation}: Define the set $\bar\calX = \{\bar x_i : i\in [1:\bar d]\}$. 

\medskip
\noindent\textit{Step 1}: \textbf{Clustering}. Choose $\calX = \{x_1,\ldots,x_d\} \subset \Omega$. For each $\ell \in [1:d]$:\\
\hspace*{1mm}$\bullet$ %
    Choose $N$~nearest neighbors $\{ x_{\ell,1},\ldots,x_{\ell,N} \} \subset \overline{\calX}$ of $x_\ell$ {with associated controls $\{u_{\ell,1},\ldots,u_{\ell,N}\}$} and define
    \begin{align}\label{eq:U_ell}
    U_\ell := 
    \begin{bmatrix}
    1 & \dots & 1\\
    u_{\ell,1} & \dots & u_{\ell,N}
    \end{bmatrix}.
    \end{align}
\hspace*{1mm}$\bullet$ %
    Approximate $g_0(x_\ell)$ and $G(x_\ell)$ by solving
    \begin{align}\label{e:argmin}
        H_\ell^*
        := \argmin\nolimits_{[y\;\;\;Y]\in\R^{n\times m+1}}\big\|\begin{bmatrix}x_{\ell,1}^+ & \dots & x_{\ell,N}^+\end{bmatrix} - \begin{bmatrix}y & Y\end{bmatrix}U_\ell \big\|_F.
    \end{align}

\noindent\textit{Step 2}: \textbf{Interpolation}. For $p\in [1:n]$ and $q\in [1:m+1]$, compute
\begin{align}\label{e:control_approx}
    \widehat H_{pq} := ((H^*_\ell)_{pq})_{\ell\in [1:d]}^\top\bK_\calX^{-1}\bk_\calX.
\end{align}
{Define $\widehat{g}_0\in \R^n$, $\widehat{G}\in \R^{n\times m}$ by %
$(\widehat{g}_0)_{p} = \widehat{H}_{p,1}$ for $p\in [1:n]$, $(\widehat{G})_{p,q} = \widehat{H}_{p,q+1}$ for $(p,q) \in [1:n]\times [1:m]$
and define $\wh f$ as in \cref{eq:dynamics:control_surr}.} %
\end{algorithm}

 We briefly comment on \cref{alg:cap} w.r.t.\ sampling and clustering: $x_1,\ldots,x_d$ can be seen as points on a ``macroscopic'' scale
    whereas the triples $(x_{\ell, i}, u_{\ell, i}, x_{\ell, i}^+)$, based on the nearest neighbors, could be understood as data on a ``microscopic level''. 
    For the sampling of the data triples $(\bar{x}_i, \bar{u}_i, \bar{x}_i^+)$, $i \in [1:\bar{d}]$, and the clustering step in \cref{alg:cap}, various strategies could be applied. 
    A canonical choice for the centers~$x_\ell$, $\ell \in [1:d]$, is to minimize the fill distance $h_\calX$, {e.g., by uniform gridding or a Padua grid as used in this work}. 
    Clearly, this choice also influences the \textit{precision} for the approximation step~\cref{e:argmin}, where we approximate the function values at $\calX$ by measurements of the flow at the nearest neighbors. %
    Both quantities, i.e., the fill distance $h_\calX$ and the radii of the clusters %
    are reflected in the %
    bound of \cref{t:control_main}. For fixed $N$, these quantities may be rendered arbitrarily small, if a sufficiently fine resolution by the micro level data set $\overline{\calX}$ is provided (e.g.~by grid-based or i.i.d.~sampling). 

Next, we present the main result of this section which provides a bound on the full approximation error $\| f-\wh f \|$. In its proof, we first estimate the approximation error made in {\it Step 1} of \cref{alg:cap} and then incorporate this bound into an error analysis concerning the interpolation in {\it Step 2}.

\begin{theorem}\label{t:control_main}
    Let {$k \geq 1+\frac{(-1)^n+1}{2}$ be the smoothness degree of the kernel and} $U_\ell\in\R^{(m+1)\times N}$ in \cref{eq:U_ell} of full rank $m+1$, $\ell\in [1:d]$. Then, there are constants $C,h_0>0$ such that for any $\calX \subset \Omega$ chosen in \cref{alg:cap} with $h_\calX < h_0$, the error between the controlled map $f$ and its surrogate $\wh f$ constructed by \cref{alg:cap} satisfies
    \begin{align*}
        \|f(x,u) - \wh f(x,u)\|_\infty \le CD(x)(1+\|u\|_1) \qquad \forall\,(x,u) \in \Omega \times \bbU,
    \end{align*}
    where 
    \begin{align*}
    D(x) &= \sqrt{2N}\big(\max\nolimits_{\ell\in [1:d]}\|U_\ell^\dagger\|\big)(L_g + L_GR)\Phi_{n,k}^{1/2}(0)\big(\max\nolimits_{v\in\mathbbm{1}}v^\top \bK_\calX^{-1}v\big)^{1/2}\varepsilon \\
    &\hspace{3cm}+ h_\calX^{k-1/2} \dist(x,\calX)\max\nolimits_{p,q}\|H_{pq}\|_{\bbH}
    \end{align*}
    with 
    $H = [g_0\ G]$, $R = \max\{\|u\|_\infty : u\in\bbU\}$, 
    maximal cluster size $\varepsilon = \max_{(\ell,j)\in [1:d]\times[1:N]} \|x_{\ell,j} - x_\ell\|$, 
    and $\mathbbm{1} = \{v\in\R^d : v_i\in\{\pm 1\}\}$.
\end{theorem}
In the proof of \cref{t:control_main}, we make use of the following lemma.
\begin{lemma}\label{l:kernel_stuff}
    For $\alpha\in\R^d$, $f = \alpha^\top\bk_\calX\in V_\calX$, $|f(x)|^2\,\le\,\k(x,x)\cdot\alpha^\top\bK_\calX\alpha
$ $\forall\, x \in \Omega$.
\end{lemma}
\begin{proof}
As $f = \sum_{k=1}^d\alpha_k\phi_{x_k}$, we have
$$
\|f\|_\bbH^2 = \langle \sum_{k=1}^d\alpha_k\phi_{x_k},\sum_{\ell=1}^d\alpha_\ell\phi_{x_\ell}\rangle = \sum_{k,\ell=1}^d\alpha_k\alpha_\ell\cdot\k(x_k,x_\ell) = \alpha^\top\bK_\calX\alpha.
$$
Consequently, the claim follows from $|f(x)|^2 = |\langle f,\phi_x\rangle_\bbH|^2\le\|\phi_x\|_\bbH^2\|f\|_\bbH^2$ and $\|\phi_x\|_\bbH^2 = \langle\phi_x,\phi_x\rangle = \phi_x(x) = \k(x,x)$.
\end{proof}

\begin{proof}[Proof of \cref{t:control_main}]
    Let $\ell \in [1:d]$. First, we show that
    \begin{align}\label{e:step1a}
        |H_{pq}(x_\ell) - (H_\ell^*)_{pq}|\,\le\,\sqrt{2N}(L_g + L_GR)\|U_\ell^\dagger\|\cdot\varepsilon_\ell, \qquad \varepsilon_\ell = \max_{j\in [1:N]} \|x_{\ell,j} - x_\ell\|.
\end{align}
For this, we have $H_\ell^* = X_\ell^+U_\ell^\dagger$, 
where $X_\ell^+ := [ x_{\ell,1}^+\ \cdots\ x_{\ell,N}^+]$. 
Hence, {since $U_{\ell}$ has full row rank, $U_\ell U_\ell^\dagger = I$ such that}
\begin{eqnarray*}
    & & |H_{pq}(x_\ell) - (H_\ell^*)_{pq}|^2 = \big|e_p^\top H(x_\ell)U_\ell U_\ell^\dagger e_q - e_p^\top X_\ell^+U_\ell^\dagger e_q\big|^2 \\
    & \le & \|e_p^\top H(x_\ell)U_\ell - e_p^\top X_\ell^+\|^2\cdot\|U_\ell^\dagger e_q\|^2 \leq \|U_\ell^\dagger\|^2 \sum\nolimits_{s=1}^{N} |e_p^\top H(x_\ell)Ue_s - e_p^\top X_\ell^+e_s |^2 \\
    & = & \|U_\ell^\dagger\|^2 \sum\nolimits_{s=1}^{N} \big| e_p^\top
        \begin{bmatrix}g_0(x_\ell)-g_0(x_{\ell,s}) & G(x_\ell) - G(x_{\ell,s})\end{bmatrix} \begin{bmatrix}1\\u_{\ell,s}\end{bmatrix} \big|^2\\
    & \le & 2\|U_\ell^\dagger\|^2\sum\nolimits_{s=1}^N \Big(\|g_0(x_\ell)-g_0(x_{\ell,s})\|^2 + \|[G(x_\ell)-G(x_{\ell,s})]u_{\ell,s}\|^2\Big)\\
    &\le & 2\|U_\ell^\dagger\|^2N(L_g^2 + L_G^2R^2)\varepsilon_\ell^2,
\end{eqnarray*}
and \cref{e:step1a} follows. Next, for $x\in\Omega$, we have
$$
    | H_{pq}(x) - \widehat H_{pq}(x) | \le | H_{pq}(x) - P_\calX H_{pq}(x) | + | P_\calX H_{pq}(x) - \widehat H_{pq}(x) |.
$$
We estimate the first summand by $\big|H_{pq}(x) - P_\calX H_{pq}(x)\big|\le Ch_\calX^{k-1/2}\dist(x,\calX)\|H_{pq}\|_{\bbH}.
$ using \cref{thm:bounds:proportional}.
For the second, we estimate
\begin{align*}
\big|P_\calX H_{pq}(x) - \widehat H_{pq}(x)\big|
&= \big|\big[(H_{pq})_\calX - [(H_\ell^*)_{pq}]_{\ell=1}^d\big]^\top\bK_\calX^{-1}\bk_\calX(x)\big|\\
&\le \big\|(H_{pq})_\calX - [(H_\ell^*)_{pq}]_{\ell=1}^d\big\|_\infty\big\|\bK_\calX^{-1}\bk_\calX(x)\big\|_1\\
&\le \sqrt{2N}(L_g + L_GR)\big(\max\nolimits_{\ell\in [1:d]}\|U_\ell^\dagger\|\big)\veps\cdot \big\|\bK_\calX^{-1}\bk_\calX(x)\big\|_1.
\end{align*}
To bound the last term, fix $x\in \Omega$ and 
$$
    \|\bK_\calX^{-1} \textbf{k}_\calX(x)\|_1 = \sum\nolimits_{i=1}^d {|e_i^\top \bK_\calX^{-1}\textbf{k}_\calX(x)|} = \sum\nolimits_{i=1}^d(-1)^{\ell_i}e_i^\top \bK_\calX^{-1}\textbf{k}_\calX(x),
$$
where $\ell_i \in \{0,1\}$ depends on~$x$. Let $v_x = \sum_{i=1}^d(-1)^{\ell_i}e_i$ and set $h := v_x^\top \bK_\calX^{-1}\textbf{k}_\calX\in V_\calX$. Then $h(x) = \|\bK_\calX^{-1}\textbf{k}_\calX(x)\|_1$ for the fixed $x$, and thus, by \cref{l:kernel_stuff}, $\|\bK_\calX^{-1}\textbf{k}_\calX(x)\|_1^2 = h^2(x) \le \k(x,x)\cdot v_x^\top \bK_\calX^{-1}v_x$. Hence, noting that for the Wendland kernels it holds that $\k(x,x) = \Phi_{n,k}(0)$, we may estimate $\|\bK_\calX^{-1}\bk_\calX(x)\|_1^2\le \Phi_{n,k}(0)\cdot\big(\max_{v\in\mathbbm{1}}v^\top \bK_\calX^{-1}v\big)$, 
which shows $\big|H_{pq}(x) - \widehat H_{pq}(x)\big|
\le C\cdot D_{pq}(x)$, where $D_{pq}(x) = c \Phi_{n,k}^{1/2}(0) ( \max\nolimits_{v\in\mathbbm{1}}v^\top \bK_\calX^{-1}v )^{1/2} \cdot \max\nolimits_{\ell\in [1:d]}\|U_\ell^\dagger\| \varepsilon$ $+$ $h_\calX^{k-1/2} \dist(x,\calX)\|H_{pq}\|_{\bbH}$ {with $c = \sqrt{2N}(L_g + L_GR)$}.
Finally, for $x\in\Omega$ and $u\in\bbU$ the claim follows from
\begin{align*}
    & \|f(x,u) - \wh f(x,u)\|_\infty
    \le \|g_0(x) - \wh g_0(x)\|_\infty + \|[G(x) - \wh G(x)]u\|_\infty\\
    \le & \max\nolimits_p|H_{p,1}(x) - \wh H_{p,1}(x)| + \max\nolimits_p \big|\sum\nolimits_{q=2}^{m+1}\big[H_{pq}(x) - \wh H_{pq}(x)]u_{q-1}\big| \\
    \leq & C \max\nolimits_{p} D_{p,1}(x) + C\max\nolimits_{p,q}D_{p,q}(x) \sum\nolimits_{q=1}^m |u_q| \leq C\max\nolimits_{p,q}D_{p,q}(x)(1+\|u\|_1).
\end{align*}
\end{proof}
{Note that the full rank condition on the controls may be viewed as a condition ensuring \emph{persistently exciting} inputs as commonly used in data-driven control, see e.g.\ \cite{de2019formulas}.}

Currently, the upper bound presented in \cref{t:control_main} depends on $h_\calX$, $N$, and~$\varepsilon$. As indicated in the discussion following \cref{alg:cap}, for fixed $N$, more data points in the micro grid $\overline{\calX}$ may be leveraged by using also a finer macro grid $\calX$ (thus, decreasing $h_\calX$) and decreasing the cluster radius $\varepsilon$, as more neighbors are contained in a smaller neighborhood. 

{If one choose choose sufficiently many control samples at each data point in \cref{alg:cap}, we may choose $\varepsilon=0$ and obtain the following corollary as a direct consequence of \cref{t:control_main}. %
\begin{corollary}\label{cor:epsnull}
Let the assumptions of \cref{t:control_main} hold, assume that we have samples of $f(x,u)$ at $(x_\ell,u_{\ell,i})$, $\ell\in [1:d]$, $i\in [1:N]$ and set $x_{\ell,i} = x_\ell$ in Step 2 of \cref{alg:cap}. Then there are $C, h_0>0$ such that for $h_{\mathcal{X}}<h_0$ and any $(x,u)\in \Omega \times \mathbb{U}$,
\begin{align*}
\|f(x,u) - \wh f(x,u)\|_\infty
\le C\cdot h_\calX^{k-1/2}\dist(x,\calX)\max\nolimits_{p,q}\|H_{pq}\|_{\bbH}(1+\|u\|_1).
\end{align*}
\end{corollary}}
We conclude this section by a remark concerning the dependence of the error bound on the number cluster elements $N$ and a possible regularization.
\begin{remark}\label{rem:thm:kEDMD-bound-control}
\noindent (a) Let us briefly discuss the behavior of the term $\sqrt{N}\|U_\ell^\dagger\|$ in the estimate in \cref{t:control_main} when the number of clusters $d$ remains constant and the cluster size $N$ grows. For this, we let $v_i := \left[\begin{smallmatrix}1\\u_{\ell,i}\end{smallmatrix}\right] \in\R^{m+1}$, $Y := \sum_{s=1}^N v_{i}v_{i}^\top$, and observe that
$U_\ell U_\ell^\top = Y$ and $U_\ell^\dagger = U_\ell^\top Y^{-1}$. Therefore $\|U_\ell^\dagger\|^2 = \|(U_\ell^\dagger)^\top U_\ell^\dagger\| = \|Y^{-1}U_\ell U_\ell^\top Y^{-1}\| = \|Y^{-1}\| = \frac 1{\lambda_{\min}(Y)}$,
where $\lambda_{\min}(Y)$ denotes the smallest eigenvalue of $Y$. Hence, we have $\sqrt N\|U_\ell^\dagger\| = \sqrt{N/\la_{\min}(Y)}$. {In what follows, we will show that if the $u_i$ are drawn independently and uniformly from $[-R,R]^m$ and $c > \max\{1,\frac {\sqrt 3}{R}\}$, then
\begin{align}\label{e:prob}
\mathbb P\big(\sqrt N\|U_\ell^\dagger\|\ge c\big)\,\le\,(m+1)\exp\left[-\tfrac{N\min\{1,\frac{R^2}{3}\}}{1+mR^2}\cdot\tfrac{\tau - \log\tau - 1}{\tau}\right],
\end{align}
where $\tau := c^2\min\{1,\frac{R^2}3\} > 1$. Hence, it is exponentially unlikely that $\sqrt N\|U_\ell^\dagger\|$ is large if $N$ is large. To prove \eqref{e:prob}, let us set $X_i = v_iv_i^\top$. Then $\|X_i\| = \|v_i\|^2 = 1 + \|u_i\|^2\le 1+mR^2 =: L$ and
$\bE[X_i] = \operatorname{diag}(1,\bE[u_iu_i^\top]) = \operatorname{diag}(1,\frac{R^2}3 I_m)$, so that $\mu_{\min} := \la_{\min}(\bE[Y]) = N\la_{\min}(\bE[X_i]) = N\min\{1,\frac{R^2}3\}$. We now have
\begin{align*}
\mathbb P\big(\sqrt N\|U_\ell^\dagger\|\ge c\big) = \mathbb P\big(\la_{\min}(Y)\le N/c^2\big) = \mathbb P\Big(\la_{\min}(Y)\le \delta\mu_{\min}\Big),
\end{align*}
where $\delta = 1/\tau < 1$. Finally, \cite[Theorem 5.1.1]{Tropp15} provides the following upper bound of the above probability, which proves \eqref{e:prob}:
\begin{align*}
(m+1)\left[e^{-(1-\delta)-\delta\log\delta}\right]^{\nicefrac{\mu_{\min}}{L}} = 
(m+1)\exp\left(-\tfrac{\mu_{\min}}L\cdot\tfrac{\frac 1\delta - \log\frac 1\delta - 1}{\frac 1\delta}\right).
\end{align*}
}

\smallskip\noindent
(b) {To alleviate a bad conditioning of the system \cref{e:control_approx} for small fill distances, one may, in addition to the clustering procedure, consider the regularized kernel matrix analogously to \cref{eq:koopman_reg} at the cost of an additional error term proportional to the regularization parameter. }

\end{remark}

\subsection{{Numerical example: Open-loop error bound}}
We illustrate the result of \cref{t:control_main} by means of a numerical example.
Therefore, consider the dynamics {obtained by an explicit Euler discretization of a controlled Duffing oscillator} given by
\begin{align}\label{eq: Duffing}
    x^+ = f(x, u) = %
    \begin{bmatrix} {x}_1 + \Delta t x_2 \\ {x}_2 + \Delta t x_1 \end{bmatrix} + \begin{bmatrix} 0\\-\Delta t3x_1^3\end{bmatrix} := g_0(x) + G(x) u
\end{align}
with $\Delta t = 0.05$ on $\Omega = [-2, 2]^2$. 
As in \Cref{subsec:numerics}, we use the coordinate functions~$\psi_i(x) = x^\top e_i$, $i \in \{1, 2\}$, as observables and for the Wendland radial basis functions, we choose the smoothness degree~$k=1$. Further, whenever we refer to randomly drawn data samples, we consider the uniform distribution on the respective set.
For the generation of the data points, {we use the Padua points introduced in \Cref{subsec:numerics} }%
to create the grid $\calX \subset \Omega$ consisting of $d$ points that will act as the set of cluster points in the macro grid. 
With randomly chosen control values $u_i\in \bbU = [-2, 2]$, the data triples~$(x_i, u_i, x_i^+)$ for $i \in [1:d]$ {satisfying the assumption of Algorithm 1} can be assembled. 

To approximate the functions $g_0$ and $G$ at the grid points in $\calX$ as presented in Step~1 of \cref{alg:cap}, for a data point~$x_\ell \in \calX$ we choose the number of neighbors~$N=25$ and then the neighbors $x_{\ell, i}$ for $i \in [1:N]$ are randomly drawn from a ball around $x_\ell$ with radius $\varepsilon = \nicefrac{1}{d}$.\footnote{{For our considered fill distances, this choice allows to compensate the error term in the bound of \Cref{{t:control_main}} as it dampens the effect of bad conditioning of the kernel matrix, i.e., $\frac{1}{d} \sqrt{{1}/{\rho_{\text{min}}(K_\mathcal{X})}} \in \{3.01, 6.16, 12.60\}$ for $d\in\{435, 1653, 6555\}$.}} For each $x_{\ell, i}$ a random $u_{\ell, i} \in \bbU$ is also drawn which yields the data points~$(x_{\ell, i}, u_{\ell, i}, x^+_{\ell, i})$ for all $\ell \in [1:d]$ needed for the regression problem. Note that $\varepsilon$ is chosen to decrease with an increasing number of cluster points. This is to make sure, that $\varepsilon$ can compensate the factor $\big(\max_{v\in\mathbbm{1}}v^\top \bK_\calX^{-1}v\big)^{1/2}$ in the error estimate in \cref{t:control_main} in view of the decreasing smallest eigenvalue of $\bK_\calX$ when decreasing the fill distance of the macro grid $\mathcal{X}$, that is, increasing $d$.
For Step~2, the approximated values of $g_0$ and $G$ on the macro grid~$\calX$ from Step~1 are used for the interpolation as stated in \cref{alg:cap} to obtain the control-affine surrogate system~$\widehat{F}$ in \eqref{eq: Duffing}. 

\begin{figure}[htb]
    \centering
    \includegraphics[width=1.\linewidth]{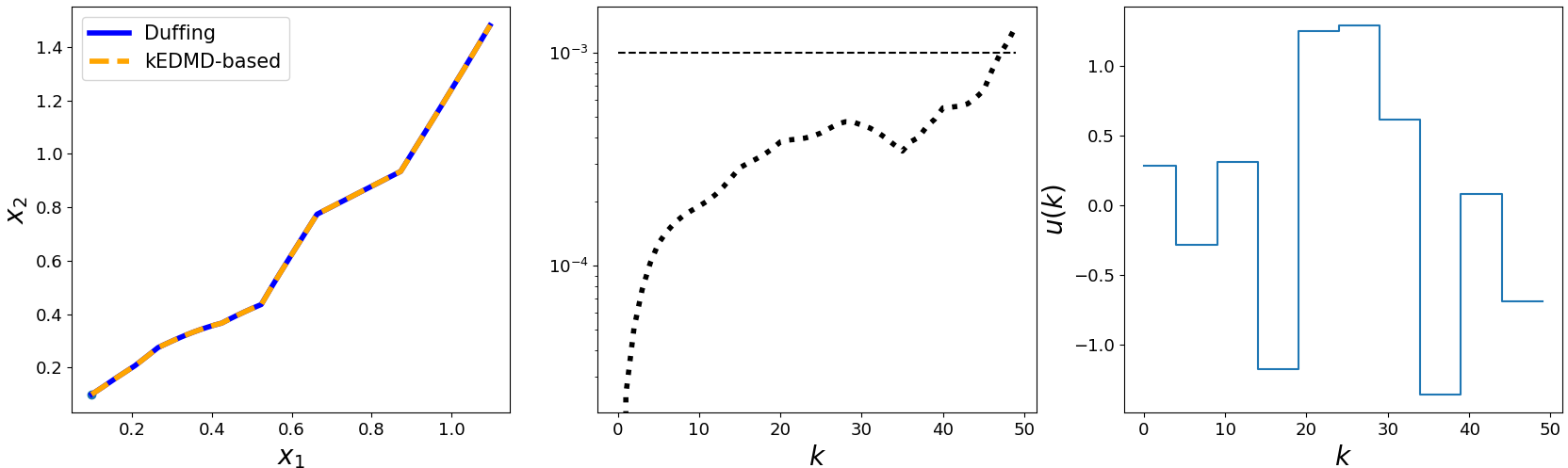}%
    \caption{Trajectories $x(\cdot;x^0,u)$ of system~\cref{eq: Duffing} and its kEDMD surrogate~$\widehat{x}(\cdot; x^0, u)$ generated with \cref{alg:cap} \braces{left} and their deviation $\| x(k; x^0, u) - \widehat{x}(k; x^0, u)\|_2$, $k \in [1:50]$, in norm {w.r.t.\ a threshold $10^{-{3}}$ }\braces{middle} for initial value $x^0 = (0.1, 0.1)^\top$ for the control sequence~$u$ \braces{right}.}
    \label{fig:kEDMD result 1}
\end{figure}

For the first simulation, we construct a kEDMD-based surrogate $\widehat{F}$ using a macro grid of $d = {6555}$ cluster points. Let $x(\cdot; x^0, u)$ and  $\widehat{x}(\cdot; x^0, u)$ denote the flow from the initial value $x^0 \in \Omega$ with control $u$ for the original model $F$ and surrogate $\widehat{F}$, respectively. In \cref{fig:kEDMD result 1}, we inspect the difference of these trajectories for the initial condition $x^0 = \binom{0.1}{0.1}$ %
and control sequences $\{u_i\}$ of length 50, where ten control values are randomly chosen from $\bbU$ and each of these control values is applied for five time steps (cf.~\cref{fig:kEDMD result 1} on the right). 
In the phase space plot in the left of \cref{fig:kEDMD result 1}, the two trajectories are barely distinguishable.
When comparing the absolute error in the middle of \cref{fig:kEDMD result 1} on a log scale, the surrogate maintains an absolute error less than {$10^{-3}$ for nearly 50 time steps} 
(corresponding to time $t=2$ in the continuous model).

Next, we randomly chose 20 initial values~$\hat{x}_i$ in $[-0.5, 0.5]^2 \subset \Omega$, as well as 20 random input sequences $\{u_i\}\subset [-0.1, 0.1]\subset\Omega$ of length 30. %
These smaller boxes ensure that the trajectories remain in the domain $\Omega$. 

\begin{figure}[htb]
    \centering
    \includegraphics[height=0.3\linewidth]{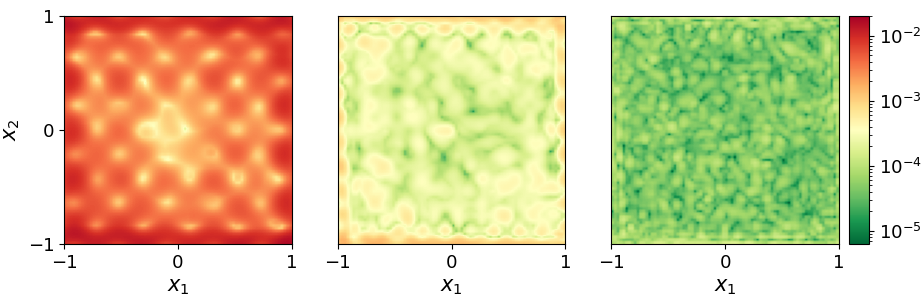}\vspace{-.3cm}
    \caption{Approximation error $\max_{j \in [1:20]}\|f(x, u_j) - \widehat{f}(x, u_j)\|$ for the predicted system~\cref{eq: Duffing} with {$d \in \{435, 1653, 6555\}$} cluster points using a {Padua}-based grid \braces{from left to right}.} %
    \label{fig:heatplot:control}
\end{figure}

Finally, in \cref{fig:heatplot:control}, we inspect the error of a one-step prediction in the Euclidean norm, similar to \cref{fig:heatplot} in \Cref{subsec:numerics}.  %
The learning of the kEDMD-based system is again performed for the three numbers of cluster points~ {$d \in \{435, 1653, 6555\}$} for the Padua-based macro grid and $N = 25$ data points in a neighborhood of $\varepsilon = \nicefrac1d$.
The errors~$\max\nolimits_{j \in [1:20]}\|F(x, u_j) - \widehat{F}(x, u_j)\|$ on a respective validation grid in $[-1, 1]^2$ are computed for the control values in $\{-2 + 0.2(j - 1) |\ j \in [1:20]\}%
$. %
In the intensity plots of \cref{fig:heatplot:control}, it can be observed that the error evenly decreases the greater $d$ is chosen.
This validates our results from \cref{t:control_main}: There are two sources of error that can be influenced, i.e., the fill distance of the macro grid~$h_\calX$ depending on number of cluster points~$d$ and the maximal cluster size~$\varepsilon$. If $\varepsilon$ is small enough to compensate the factor $\big(\max\nolimits_{v\in\mathbbm{1}}v^\top \bK_\calX^{-1}v\big)^{\nicefrac12}$ that increases with $d$, the decrease of the fill distance by refining $\calX$ yields a decay of the approximation error. 

\subsection{Feedback stabilization}\label{subsec:feedback}
Finally, we consider stabilizing feedback laws for the original control-affine system \cref{eq:dynamics:control} and its surrogate \cref{eq:dynamics:control_surr}, where $\wh f$ is defined via \cref{alg:cap}. As it turns out, under certain assumptions, a stabilizing feedback for \cref{eq:dynamics:control_surr} is also stabilizing for \cref{eq:dynamics:control}.

\begin{proposition}\label{prop:feedback}
    Let $H = [g_0\;\;G]$ {be componentwise in $C_b^{\ceil{\sigma_{n,k}}}(\Omega)$ with $\sigma_{n,k} = (n+1)/2+k$ and $k\geq 1$.} %
    Further, let $\kappa :\Omega \to \R^n$ be a {locally Lipschitz-continuous } feedback law asymptotically stabilizing %
    the equilibrium $x^* \in \Omega \cap \calX$ of the closed-loop dynamics $\widehat{F}_\mathrm{cl}(x) := \widehat f(x,\kappa(x))$ such that the modulus of continuity~$\omega_V$ of the respective Lyapunov function~$V$ and %
    is uniform in $h_\calX\in (0,h_0]$ for some $h_0>0$.  %
    Then, $x^*$ is practically asymptotically stable w.r.t.~$f(\cdot,\kappa(\cdot)) =: F_{\rm cl}$ with practical region decreasing in the fill distance $h_\calX$ of the macro grid $\calX$ and the cluster size~$\varepsilon$.
    
    If further the samples are drawn %
    with cluster size $\varepsilon=0$
    and if $\widehat{F}_\mathrm{cl}(x) := \widehat f(x,\kappa(x))$ is asymptotically stable towards $x^*\in \Omega$ with a Lyapunov function satisfying Assumption \braces{a} or \braces{b} of \cref{thm:lyapunov:kEDMD:stability}, then $\kappa$ also asymptotically stabilizes the original dynamics defined in \cref{eq:dynamics:control} towards $x^*$.
\end{proposition}
\begin{proof}
    We mimic the proofs of \cref{thm:lyapunov:kEDMD:practical} and \cref{thm:lyapunov:kEDMD:stability} for the closed-loop system. %
    In the first case, we observe that the error bound of \cref{t:control_main} yields 
    \begin{align*}
        \|F_\mathrm{cl}(x) - \widehat{F}_\mathrm{cl}(x)\| &= \|f(x,\kappa(x)) - \widehat f(x,\kappa(x))\|\\
        &\leq (c_1h_\calX^{k-1/2} \dist(x,\calX)+ c_2\varepsilon )(1+\|\kappa(x)\|_1) \leq c_3(h_\calX^{k+1/2}+ \varepsilon)
    \end{align*}
    where we invoked the boundedness of $\kappa$ and with constants $c_1 := C\max_{p,q}\|H_{pq}\|_\bbH$, $c_2 := C\sqrt{2N}(L_g + L_GR)\Phi_{n,k}^{1/2}(0)\big(\max_{\ell\in [1:d]}\|U_\ell^\dagger\|\big)\big(\max_{v\in\mathbbm{1}}v^\top \bK_\calX^{-1}v\big)^{1/2}$ and $c_3 := \max\{c_1,c_2\}(1+\sup_{x\in\Omega} \|\kappa(x)\|_1)$. This error bound is structurally the same as the bound of \cref{cor:Fbound_nonproportional} with the cluster size $\varepsilon$ taking the role of $\sqrt{\lambda}$ and $\omega_\Upsilon = \operatorname{id}$. {Thus, we obtain a similar estimate to \eqref{eq:criticalestimate} with $\bar{C}h_\mathcal{X}^{1/2}$ replaced by $c_3(h_{\mathcal{X}}^{k+1/2}+\varepsilon)$. To mimic remainder of the proof of \cref{thm:lyapunov:kEDMD:practical}, we remain to show the Lipschitz continuity of $\widehat{F}_\mathrm{cl}(x) =  \widehat{g}_0 + \widehat{G}(x)\kappa(x)$, uniformly in the fill distance $h_\mathcal{X}$. With the same argumentation as in \cref{lem:Lipschitz:uniform}, we may first show that $\widehat{g}_0$ and $\widehat{G}$ are Lipschitz continuous. As the feedback is assumed to be locally Lipschitz, and as all functions are bounded due to continuity, $\widehat{F}_\mathrm{cl}$ is also locally Lipschitz continuous as sum and product of bounded and locally Lipschitz continuous functions.} %
    Thus, a similar argumentation as in the {second part of the} proof of \cref{thm:lyapunov:kEDMD:practical} yields the claim.
    
    To prove the second claim, we note that in view of the sampling of \cref{cor:epsnull}, we may set $\varepsilon= 0$ and obtain
    \begin{align*}
        \|F_\mathrm{cl}(x) - \widehat{F}_\mathrm{cl}(x)\| %
        &\leq h_\calX^{k-1/2} \|x-x^*\| \max\nolimits_{p,q}\|H_{pq}\|_\bbH (1+\sup\nolimits_{x\in \Omega}\|\kappa(x)\|_1),
    \end{align*}
    where we used that $\dist(x,\calX) \leq \|x-x^*\|$ due to $x^*\in \calX$. This inequality in particular implies that $x^*$ is also an equilibrium of ${F}_\mathrm{cl}(x)$. Further, the above is a proportional bound as in \cref{c:surrogate_proportional}. Thus, an inspection of the proof of \cref{thm:lyapunov:kEDMD:stability} shows that this implies also asymptotic stability of $F_\mathrm{cl}$.
\end{proof}

{Exemplarily, we consider the discrete-time system from~\cite[Example~3.1]{kazantzis2001nonlinear}, i.e.,
\begin{align}\label{eq:kanzantzis}
    \begin{bmatrix}
        x_1^+ \\ x_2^+
    \end{bmatrix} = \tilde{f}(x, u) = \begin{bmatrix}
        \exp(0.3x_2)\sqrt{1 + x_1 + x_2} - 0.4x_2 - 1 \\
        0.5 \ln(1 + x_1 + x_2) + 0.4x_2
    \end{bmatrix} + \begin{bmatrix}
        0.5 \\ 0
    \end{bmatrix} u
\end{align}
with map~$\tilde{f}: \tilde\Omega %
\rightarrow \R^2$ and domain $\tilde\Omega = [-0.4, 1]^2$. For vanishing input~$u = 0$, the origin $x^*=0$ clearly is an equilibrium of system~\eqref{eq:kanzantzis}. However, as shown in~\cite{kazantzis2001nonlinear}, this equilibrium is not asymptotically  stable, as the linearization at the origin has an eigenvalue of modulus larger than one. To stabilize the origin, \cite{kazantzis2001nonlinear} introduces the stabilizing state feedback $\kappa(x) = -\ln(1 + x_1 + x_2)$.
Using this example, we illustrate that the converse statement of Proposition~\ref{prop:feedback}, i.e., if the original system~\eqref{eq:kanzantzis} is asymptotically stable towards the origin for the feedback law~$\kappa$, the origin is practically asymptotically stable w.r.t.\ the kEDMD-based surrogate system. }

\begin{figure}[htb]
    \centering
    \includegraphics[width=.96\linewidth]{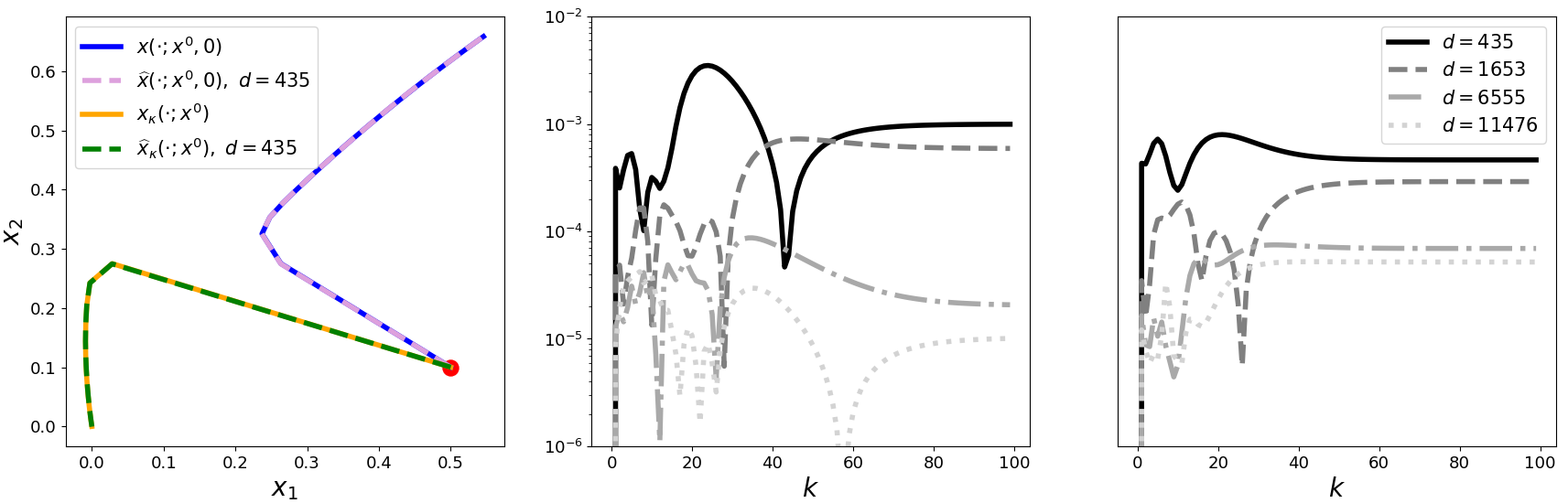}\\[-2mm]
    \caption{
    {Left: KEDMD approximants~$\widehat{x}(\cdot; x^0, 0)$ (open loop) and~$\widehat{x}_\kappa(\cdot; x^0)$ (closed loop) of system~\cref{eq:kanzantzis}, initial value $x^0 = (0.5, 0.1)^\top$, generated with \cref{alg:cap} based on the Padua grid with $d = 435$ cluster points. Deviations $\| x(k; x^0, 0) - \widehat{x}(k; x^0, 0)\|_2$ and $\| x_\kappa(k; x^0) - \widehat{x}_\kappa(k; x^0)\|_2$, $k \in [1:100]$, in open-(mid) and closed-loop (right) control for $d \in \{435, 1653, 6555, 11476\}$.}
    }
    \label{fig:Padua:control:feedback}
\end{figure}

{To this end, we approximate the system~\eqref{eq:kanzantzis} using \Cref{alg:cap}. The set of virtual observation points~$\mathcal{X}$ is chosen as a grid of $d \in \{435, 1653, 6555, 11476\}$ Padua points, as introduced in \Cref{subsec:numerics} and the cluster radius~$\varepsilon$ is again set to $\tfrac{1}{d}$. 
Then, for each point in $\mathcal{X}$, we randomly draw 25 data points~$(x_{\ell, i}, u_{\ell, i}, x_{\ell, i}^+)$ for $i \in [1:d]$ and $\ell \in [1:25]$ such that the assumptions of \Cref{alg:cap} are satisfied.  
As before, the coordinate functions are chosen as observables and the Wendland kernel with smoothness degree~$k = 1$ is used. The results are depicted in \Cref{fig:Padua:control:feedback}. For reasons of clarity, the trajectories are only shown for $d = 435$ (left), as the phase portrait does not show noticeable differences for higher numbers of virtual observation points. The errors (mid, right), however, are presented for all $d \in \{435, 1653, 6555, 11476\}$.
}

\noindent {We first inspect the open-loop predictions for the uncontrolled case. As can be seen in \Cref{fig:Padua:control:feedback} (left), starting at the initial state~$x^0 = (0.5, 0.1)^\top$ the behavior of the ground truth is accurately replicated, as both trajectories move towards the second equilibrium at $x^* = (0.55, 0.66)^\top$ illustrating that this initial value is not contained in a potential domain of attraction. 
The kEDMD-based surrogate can only approach $x^*$ up to a certain accuracy, as shown in the error plots (\Cref{fig:Padua:control:feedback}, mid). It can be observed, that the error decreases the greater the number~$d$ of Padua grid points is. 
Furthermore, \Cref{fig:Padua:control:feedback} exhibits the results for using the nominally stabilizing feedback~$\kappa$. 
The feedback law also stabilizes the origin w.r.t.\ the kEDMD approximants as can be inferred from the numerically computed trajectories of system~\cref{eq:kanzantzis}, which tend towards the origin. 
When using $\kappa$ also for the kEDMD-based surrogate systems, we observe the statement of \Cref{prop:feedback} (with switched roles of surrogate and ground truth). The trajectories approach a neighborhood of the origin depending on the number of cluster points. These results confirm that asymptotically stabilizing feedbacks for the original system are practically asymptotically stabilizing w.r.t.\ the kEDMD-based surrogate. Future work considers leveraging the converse statement (as proven in \Cref{prop:feedback}) for data-driven computation of feedback laws. 
}

\section{Conclusions and outlook}
\label{sec:conclusions}

We provided a novel kernel EDMD scheme with flexible state-control sampling and uniform error estimates for data-driven modeling of dynamical (control) systems with stability guarantees.
In the first part of this work, we extended existing uniform bounds~\cite{KohnPhil24} on the full approximation error on Koopman approximants generated with kernel EDMD to regularized ones and derived proportional error bounds. 
While the first is important for robustness in view of the poor conditioning of kernel matrices (e.g.~for noisy data), the second is key to show that asymptotic stability of an equilibrium w.r.t.\ the surrogate dynamics, certified by a Lyapunov function, is preserved for the original dynamical system and vice versa. 
In the second part, we proposed a novel kernel-EDMD scheme for control-affine systems building upon arbitrarily-sampled control-state pairs and rigorously showed the first uniform (proportional) bounds on the full approximation error for control systems. 
We demonstrated that, as a consequence, stabilizing feedback laws designed for the data-driven surrogate models also ensure asymptotic stability w.r.t.\ the closed-loop of the original system building upon the proposed stability-analysis framework for kernel EDMD. 
We accompanied our findings by various numerical simulations.

Future work will be devoted to data-driven controller design leveraging the proposed highly-flexible sampling regime in combination with the novel proportional bounds derived in a non-restrictive setting. 
Herein, results from~\cite{huang2023robust} and~\cite{dorfler2022bridging} may be leveraged for robustification and direct data-driven controller design.
{To alleviate the curse of dimensionality, %
future work considers leveraging low-dimensional structures by means of, e.g., tensor-train approximations~\cite{holtz2012alternating}, greedy methods for feature selection~\cite{santin2024optimality} or neural networks for dictionary learning~\cite{jin2024extended}.}

\bibliographystyle{plain}
\bibliography{references}

\appendix

\section{Results for the alternative surrogate}\label{s:K_tilde}
In this part of the Appendix, we discuss the properties of the approximant $\widetilde\calK = P_\calX\calK P_\calX
$ of the Koopman operator (cf.\ \Cref{rem:alternative}), where $\calX = \{x_1,\ldots,x_d\}\subset\Omega$, and $P_\calX$ is the orthogonal projection onto $V_\calX = \linspan\{\phi_{x_1},\ldots,\phi_{x_d}\}$. More generally, for $\la\ge 0$, we consider the regularized approximant
$$
\wt\calK_\la := R_\calX^\la\calK R_\calX^\la.
$$
As for most of \Cref{sec:Koopman}, we assume that the RKHS $\bbH$ is invariant under the Koopman operator $\calK$, i.e. $\calK\bbH\subset\bbH$. First of all, the approximant $\wt\calK_\la$ can be written as
$$
\widetilde\calK_\la f = f_\calX^\top(\bK_\calX + \la I)^{-1}\bK_{\calX,F(\calX)}(\bK_\calX + \la I)^{-1}\bk_\calX.
$$
Since
$$
F(x) = \Upsilon(\Psi(F(x))) = \Upsilon\big(\big[\calK\psi_1(x),\ldots,\calK\psi_M(x)\big]^\top\big),
$$
we define another family of surrogate models for \Cref{eq:dynamics} by
\begin{equation}\label{e:surrogate_hat}
x^+ = \widetilde F_\la(x),
\end{equation}
where
\begin{align*}
\widetilde F_\la(x)
&= \Upsilon\big(\big[\wt\calK_\la\psi_1(x),\ldots,\wt\calK_\la\psi_M(x)\big]^\top\big)\\
&= \Upsilon\big(\Psi_\calX^\top(\bK_\calX + \la I)^{-1}\bK_{\calX,F(\calX)}(\bK_\calX + \la I)^{-1}\bk_\calX(x)\big).
\end{align*}
In what follows, we will prove analogues of the basic results in \Cref{sec:Koopman} and \Cref{sec:stability} involving the approximant $\wh\calK_\la$ and the associated surrogate \Cref{eq:dynamics:surrogate} for the alternative approximant $\wt\calK_\la$ and the surrogate \Cref{e:surrogate_hat}, introduced above. It is then clear that our stability results in \Cref{sec:stability} hold with $\wh F_\la$ replaced by $\wt F_\la$.

\subsection{Error bounds}
The following theorem is the analogue of the combination of \Cref{thm:error:regularization} and \Cref{cor:Fbound_nonproportional} for the alternative approximants.

\begin{theorem}\label{thm:error:regularization:alt}
	Let $k\geq 1$, $\lambda \geq 0$ and $F\in C_b^{\ceil{\sigma_{n,k}}}(\Omega;\R^{n})$. Then there are constants $C,h_0>0$ such that for any finite set $\calX = \{x_i\}_{i=1}^d \subset\Omega$ of sample points with $h_{\calX}\le h_0$ and for all $\lambda \geq 0$, $f\in\calN_{\Phi_{n,k}}(\Omega)$ and $x\in\Omega$ we have
	\begin{align*}    
	|(\mathcal{K}f)(x) - (\wt{\mathcal{K}}_\lambda f)(x)| \leq C\big(h_\calX^{k+1/2} + \sqrt\lambda\big)\|f\|_{\calN_{\Phi_{n,k}}(\Omega)}.
	\end{align*}
	In particular,
	\begin{align*}    
	\|F(x) - \wt{F}_\lambda (x)\| \leq \omega_\Upsilon \Big(C\big(h_\calX^{k+1/2} + \sqrt\lambda\big)\|\Psi\|_{\calN^l_{\Phi_{n,k}}(\Omega)}\Big) \qquad \forall x\in \Omega
	\end{align*}
\end{theorem}
\begin{proof}
	The first claim on the error between $\calK$ and $\wt\calK_\la$ follows by a combination of \cite[Theorem 5.2]{KohnPhil24} and \Cref{thm:Wend:regularization}. The proof of the second claim on $F-\wt F_\la$ is very similar to that of \Cref{cor:Fbound_nonproportional}.
\end{proof}

The next theorem improves the statement of \Cref{thm:error:regularization:alt} towards proportional bounds, similarly as \Cref{thm:bounds:proportional} and \Cref{c:surrogate_proportional} did for the approximant $\wh\calK$.

\begin{theorem}[Proportional error bounds]\label{thm:bounds:proportional:alt}
	Assume that $k \geq 1+\frac{(-1)^n+1}{2}$ and $F \in C_b^{\ceil{\sigma_{n,k}}}(\Omega;\R^{n})$. Then there exist $C, h_0>0$ such that, for any finite set $\calX = \{x_i\}_{i=1}^d \subset \Omega$ of sample points with $h_{\calX}\le h_0$ and all $f \in \calN_{\Phi_{n,k}}(\Omega)$ and $x\in\Omega$ we have
	\begin{align*}    
	|(\mathcal{K}f)(x) - (\widetilde{\mathcal{K}}f)(x)| \leq Ch_\calX^{k-1/2} \max\{\dist(x,\calX),\dist(F(x),\calX)\}\|f\|_{\calN_{\Phi_{n,k}}(\Omega)}.
	\end{align*}
	In particular,
	\begin{align*}    
	\|F(x) - \wt F(x)\|\le \omega_\Upsilon\big(Ch_\calX^{k-1/2}\|\Psi\|_{\calN_{\Phi_{n,k}}(\Omega)^M}\max\{\dist(x,\calX),\,\dist(F(x),\calX)\}\big).
	\end{align*}
	Moreover, if $x^*\in \calX$ is an equilibrium of the dynamics \Cref{eq:dynamics}, then
	\begin{align*}    
	\|F(x) - \wt F(x)\|\le \omega_\Upsilon\big(Ch_\calX^{k-1/2}\|\Psi\|_{\calN_{\Phi_{n,k}}(\Omega)^M}\|x-x^*\|\big).
	\end{align*}
\end{theorem}
\begin{proof}
	We may use \Cref{thm:bounds:proportional} to estimate
	\begin{align*}
	|(\mathcal{K} f)(x) - (\widetilde{\mathcal{K}} f)(x)| & \le |(\mathcal{K} f)(x) - (\mathcal{K} P_\calX f)(x)| + |(\mathcal{K} P_\calX f)(x) - (P_\calX\mathcal{K} P_\calX f)(x)| \\
	&= |[\mathcal{K}(I-P_\calX)f](x)| + |[(I-P_\calX)\mathcal{K} P_\calX f](x)|\\
	& \le |[(I - P_\calX)f](F(x))| + C h_{\calX}^{k-1/2}\dist(x,\calX) \|\mathcal{K} P_\calX f\|_{\calN_{\Phi_{n,k}}} \\
	&\le C h_{\calX}^{k-1/2}\!\Big[\!\dist(F(x),\calX) \!+\! \dist(x,\calX) \|\mathcal{K}\|_{\calN_{\Phi_{n,k}}\to\calN_{\Phi_{n,k}}}\!\Big]\|f\|_{\calN_{\Phi_{n,k}}},
	\end{align*}
	which proves the main statement. The first claim of the ``in particular''-part is now clear, the second follows from $\dist(x,\calX)\le \|x-x^*\|$ and
	$$
	\dist(F(x),\calX)\le \|F(x)-x^*\| = \|F(x)-F(x^*)\|\le \|F\|_{C_b^1(\Omega,\R^{n})}\|x-x^*\|.
	$$
	Hence, the theorem is proved.
\end{proof}

\subsection{Equilibria}
Let us compare the equilibria in the data set $\calX$ of the dynamics \Cref{eq:dynamics} and the surrogate \Cref{e:surrogate_hat}.

\begin{lemma}\label{lem:equi}
	If $x^*\in\calX$ is an equilibrium of \Cref{eq:dynamics}, then for each $f\in\bbH$ we have
	\[
	(\widetilde{\mathcal{K}}f)(x^*) = (\mathcal{K}f)(x^*) = f(x^*).
	\]
	In particular, $x^*$ is an equilibrium of \Cref{e:surrogate_hat}.
\end{lemma}
\begin{proof}
	Let $x^* = x_k$, and let $f\in\bbH$ be arbitrary. Clearly, $(\mathcal{K} f)(x_k) = f(F(x_k)) = f(x_k)$. Note that $P_\calX h$ for $h\in\bbH$ is the unique function in $V_\calX$ which satisfies $(P_\calX h)(x) = h(x)$ for all $x\in\calX$. Hence, if we set $g = \mathcal{K} P_\calX f$, then $(P_\calX g)(x_k) = g(x_k)$ as $g\in\bbH$ and $x_k\in\calX$. Therefore,
	\[
	(\widetilde{\mathcal{K}} f)(x_k) \!=\! (P_\calX g)(x_k) \!=\! g(x_k) \!=\! (\mathcal{K} P_\calX f)(x_k) \!=\! (P_\calX f)(F(x_k)) \!=\! (P_\calX f)(x_k) = f(x_k),
	\]
	which proves the lemma.
\end{proof}

The next proposition characterizes the equilibria in the data set $\calX$ of \Cref{eq:dynamics} and the surrogate \Cref{e:surrogate_hat}.

\begin{proposition}
	Let $x^* = x_k\in\calX$. Then the following hold:
	\begin{enumerate}
		\item[{\rm (i)}]  If $\Psi_\calX^\top(I - \bK_\calX^{-1}\bK_{\calX,F(\calX)})e_k=0$, then $x^*$ is an equilibrium of \Cref{e:surrogate_hat}. The converse is true if $\Upsilon : \R^M\to\R^n$ is injective.
		\item[{\rm (ii)}] If $x^*$ is an equilibrium of \Cref{eq:dynamics}, then $(I - \bK_\calX^{-1}\bK_{\calX,F(\calX)})e_k=0$. The converse holds for kernels $\k$ with constant diagonal, in particular for the Wendland kernels.
	\end{enumerate}
\end{proposition}
\begin{proof}
	(i). We have
	$$
	\widetilde F(x_k) = \Upsilon\big(\Psi_\calX^\top\bK_\calX^{-1}\bK_{\calX,F(\calX)}\bK_\calX^{-1}\bk_\calX(x_k)\big) = \Upsilon\big(\Psi_\calX^\top\bK_\calX^{-1}\bK_{\calX,F(\calX)}e_k\big)
	$$
	Hence, if $\Psi_\calX^\top(I - \bK_\calX^{-1}\bK_{\calX,F(\calX)})e_k=0$, then
	$$
	\widetilde F(x_k) = \Upsilon\big(\Psi_\calX^\top e_k\big) = \Upsilon(\Psi(x_k)) = x_k.
	$$
	Conversely, if $\widetilde F(x_k) = x_k$, then
	$$
	\Upsilon\big(\Psi_\calX^\top\bK_\calX^{-1}\bK_{\calX,F(\calX)}e_k\big) = x_k = \Upsilon(\Psi(x_k)).
	$$
	Hence, if $\Upsilon$ is injective, then $\Psi_\calX^\top\bK_\calX^{-1}\bK_{\calX,F(\calX)}e_k = \Psi(x_k) = \Psi_\calX^\top e_k$, which implies the claim.
	
	\smallskip\noindent
	(ii). Let $x^*=x_k$ be an equilibrium of \Cref{eq:dynamics}. Then
	$$
	(I - \bK_\calX^{-1}\bK_{\calX,F(\calX)})e_k = e_k - \bK_\calX^{-1}\bk_\calX(F(x_k)) = e_k - \bK_\calX^{-1}\bk_\calX(x_k) = e_k - \bK_\calX^{-1}\bK_\calX e_k = 0.
	$$
	Conversely, let $(I - \bK_\calX^{-1}\bK_{\calX,F(\calX)})e_k = 0$, and assume that $k(x,x) = c > 0$ for all $x\in\Omega$. Then $(\bK_\calX - \bK_{\calX,F(\calX)})e_k = 0$, i.e.,
	$$
	\bk_\calX(x^*) = \bk_\calX(y^*),
	$$
	where $y^* = F(x^*)$. In particular, $\k(x^*,y^*) = \k(x^*,x^*) = c$. But this means that $|\langle \phi_{x^*},\phi_{y^*}\rangle|^2 = c^2 = \k(x^*,x^*)\k(y^*,y^*) = \|\phi_{x^*}\|^2\|\phi_{y^*}\|^2$. By Cauchy-Schwarz, $\phi_{x^*}$ and $\phi_{y^*}$ must be linearly dependent, which (by the strict positive definiteness of the kernel $\k$) implies $y^* = x^*$.
\end{proof}

\section{The regularization operator}
Recall that the linear operator $R_\calX^\lambda$ on $\bbH$ is defined by
$$
R_\calX^\lambda f = f_\calX^\top(\bK_\calX + \lambda I)^{-1}\bk_\calX.
$$
The following proposition is not used in the paper, but might be of independent interest.

\begin{proposition}\label{p:RX}
	The following statements hold for the operator $R_\calX^\lambda : \bbH\to\bbH$:
	\begin{enumerate}
		\item[(i)]   $(R_\calX^\lambda)^* = R_\calX^\lambda\,\ge\,0$
		\item[(ii)]  $R_\calX^\lambda P_\calX = P_\calX R_\calX^\lambda = R_\calX^\lambda$
		\item[(iii)] $\|R_\calX^\lambda f\|\le\|P_\calX f\|\le\|f\|$ for $f\in\bbH$.
	\end{enumerate}
\end{proposition}
\begin{proof}
	(i). Let $f,g\in\bbH$. Then, with $\alpha = (\bK_\calX+\lambda I)^{-1}f_\calX$ we have
	\begin{align}\label{e:nett}
	\langle R_\calX^\lambda f,g\rangle = \Big\langle\sum_{i=1}^n\alpha_i\phi(x_i),g\Big\rangle = \sum_{i=1}^n\alpha_ig(x_i) = \alpha^\top g_\calX = f_\calX^\top (\bK_\calX+\lambda I)^{-1}g_\calX.
	\end{align}
	Similarly, one shows that $\langle f,R_\calX^\lambda g\rangle = f_\calX^\top(\bK_\calX + \lambda I)^{-1}g_\calX$. This shows $R_\calX^\lambda = (R_\calX^\lambda)^*\,\ge\,0$.
	
	\medskip
	(ii). $P_\calX R_\calX^\lambda = R_\calX^\lambda$ follows from $\operatorname{ran} R_\calX^\lambda\subset V_\calX$ and $R_\calX^\lambda P_\calX = R_\calX^\lambda$ from $R_\calX^\lambda|_{V_\calX^\perp} = 0$. For the latter, note that $V_\calX^\perp = \{f\in\bbH : f(x_k)=0\,\forall k\in [1:d]\}$.
	
	\medskip
	(iii). Let $g=R_\calX^\lambda f$. Then $g_\calX = \bK_\calX(\bK_\calX+\lambda)^{-1}f_\calX$, and thus plugging $g$ into \Cref{e:nett} yields
	\begin{align*}
	\|R_\calX^\lambda f\|^2
	&= f_\calX^\top(\bK_\calX+\lambda)^{-1}\bK_\calX(\bK_\calX+\lambda)^{-1}f_\calX\\
	&= f_\calX^\top\big[\underbrace{(\bK_\calX+\lambda)^{-1}\bK_\calX(\bK_\calX+\lambda I)^{-1} - \bK_\calX^{-1}}_{=:\mathbf{A}}\big]f_\calX + f_\calX^\top\bK_\calX^{-1}f_\calX\\
	&= f_\calX^\top\mathbf{A} f_\calX + \|P_\calX f\|^2.
	\end{align*}
	By the spectral mapping theorem, the set of eigenvalues of $\mathbf{A}$ equals
	\[
	\sigma(\mathbf{A}) = \left\{\frac{\mu}{(\mu+\lambda)^2} - \frac 1{\mu} : \mu\in\sigma(\bK_\calX)\right\} = \left\{\frac{\mu^2-(\mu+\lambda)^2}{\mu(\mu+\lambda)^2} : \mu\in\sigma(\bK_\calX)\right\}.
	\]
	Hence, all eigenvalues of $\mathbf{A}$ are negative, which implies that $\mathbf{A}$ is negative definite. Therefore, $\|R_\calX^\lambda f\|^2 = f_\calX^\top\mathbf{A}f_\calX + \|P_\calX f\|^2\le\|P_\calX f\|^2$. The second inequality $\|P_\calX f\|\le\|f\|$ is clear as $P_\calX$ is an orthogonal projection.
\end{proof}

\end{document}